\documentclass[twocolumn]{autart} 
\usepackage{graphicx} 
\usepackage{amssymb}
\usepackage{amsmath}
\usepackage{subfigure}

\begin{document}

\begin{frontmatter}

\title{Solving A Class of Nonsmooth Resource Allocation Problems with Directed Graphs though Distributed Smooth Multi-Proximal Algorithms}

\thanks[footnoteinfo]{This work was supported by Projects of Major International (Regional) Joint Research Program NSFC (Grant no. 61720106011), NSFC (Grant no.61621063, 61573062, 61673058, 61873033), Beijing Advanced Innovation Center for Intelligent Robots and Systems (Beijing Institute of Technology), Key Laboratory of Biomimetic Robots and Systems (Beijing Institute of Technology), Ministry of Education, Beijing, 100081, China. P. Pardalos was supported by the Paul and Heidi Brown Preeminent Professorship at ISE, University of Florida and a Humboldt Research Award. Corresponding author: Hao Fang.}

\author[PC,BIT,QH]{Yue Wei}\ead{weiy@pcl.ac.cn},    
\author[BIT]{Chengsi Shang} \ead{chengsi.shang@qq.com},
\author[BIT]{Hao Fang}\ead{fangh@bit.edu.cn},              
\author[BIT]{Xianlin Zeng}\ead{xianlin.zeng@bit.edu.cn}, 
\author[BIT]{Lihua Dou}\ead{doulihua@bit.edu.cn}, 
\author[UOF]{Panos Pardalos}
\ead{pardalos@ufl.edu}

\address[PC]{Peng Cheng Laboratory, China}
\address[BIT]{The Key Laboratory of Intelligent Control and Decision of Complex Systems, Beijing Institute of Technology, China}  %
\address[QH]{Graduate School at Shenzhen, Tsinghua University, China}
\address[UOF]{Center of Applied Optimization (CAO), Industrial and Systems Engineering, University of Florida, USA}


\begin{abstract}
In this paper, two distributed multi-proximal primal-dual algorithms are proposed to deal with a class of distributed nonsmooth resource allocation problems. In these problems, the global cost function is the summation of local convex and nonsmooth cost functions, each of which consists of one twice differentiable function and multiple nonsmooth functions. Communication graphs of underling multi-agent systems are directed and strongly connected but not necessarily weighted-balanced. The multi-proximal splitting is designed to deal with the difficulty caused by the unproximable property of the summation of those nonsmooth functions. Moreover, it can also guarantee the smoothness of proposed algorithms. Auxiliary variables in the multi-proximal splitting are introduced to estimate subgradients of nonsmooth functions. Theoretically, the convergence analysis is conducted by employing Lyapunov stability theory and integral input-to-state stability (iISS) theory with respect to set. It shows that proposed algorithms can make states converge to the optimal point that satisfies resource allocation conditions.
\end{abstract}

\end{frontmatter}

\section{Introduction}

In this paper, we consider a class of distributed nonsmooth convex resource allocation problems with directed graphs. A wide range of problems in the field of coordination of multi-agent systems \cite{ZY2017}-\cite{JC2018}, economic dispatch of power systems \cite{Cortes2016} and machine learning belong to this class of problems. As examples, in distributed constrained coordination of multi-agent systems with directed graphs, the local cost function of agent $i$ usually consists of a smooth function and multiple nonsmooth functions standing for different constraints and tasks. Moreover, multi-agent systems are required to maintain some configurations described by resource allocation conditions. When considering a classical machine learning problem - the fused LASSO problem \cite{Lasso} - with constraints and directed graphs, the least squares loss is smooth. The $l_{1}$ penalty and indicator functions of local constraints in this problem are usually nonsmooth. Then resource allocation conditions are employed here as global constraints. As common features, each global cost function in these problems is summed up by local cost functions, and each local cost function consists of a smooth convex function and multiple nonsmooth convex functions. Even though nonsmooth functions are proximable, their summation might not be, where a function being proximable means that the proximal operator of this function has a closed or semi-closed form solution and is computationally easy to evaluate \cite{CTOS3}. Besides, connected graphs of these problems are directed and maybe weight-unbalanced, where a directed graph being weight-unbalanced means that the in-degree and out-degree of some nodes in this graph are unequal. The difficulty of these problems is to tackle nonsmooth cost functions and directed connecting graphs simultaneously. Due to important applications and challenges mentioned above, these problems have attracted increasing attentions.

\subsection*{Literature review}

Communication between agents in multi-agent
systems has attracted much attention due to the importance of information exchange. Recently, continuous-time distributed algorithms for resource allocation problems have been widely investigated with different kinds of connected graphs \cite{YP2016}-\cite{YWW}. For undirected graphs, \cite{YP2016} designed an initialization-free distributed algorithm for distributed resource allocation problems. \cite{YH2019} proposed a new distributed private-guaranteed algorithm to solve economic dispatch problems with undirected graphs. As to directed graphs, \cite{SP} proposed a continuous-time algorithm via singular perturbation for distributed resource allocation problems. While \cite{SP} did not consider local constraints. In \cite{SP1}, a distributed projection-based algorithm was designed to deal with distributed resource allocation problems with weight-balanced graphs. \cite{DZ2019} investigated constrained nonsmooth resource allocation problems via a distributed algorithm, which can solve resource allocation problems with strongly convex cost functions and weight-balanced digraphs, as well as resource allocation problems with strictly convex cost functions and connected undirected graphs. For distributed resource allocation problems with weight-unbalanced graphs, \cite{YWW} proposed a distributed adaptive algorithm to achieve the optimal solution. While this algorithm fails to solve resource allocation problems with local constraints and weight-unbalanced graphs simultaneously. 

Nonsmoothness is a natural property of many resource allocation problems in real-world science and engineering areas. Two important categories of existing algorithms for solving distributed nonsmooth optimization and resource allocation problems are shown here. The first category is subgradient-based algorithms proposed in \cite{DCO1}-\cite{DCO5}, whose convergence was proven based on nonsmooth analysis \cite{NCO}. \cite{Xie1} designed a distributed continuous-time projected algorithm to deal with distributed constrained nonsmooth optimization problems. \cite{DOC1} investigated the distributed nonsmooth constrained optimization problem with distributed projection-based saddle-point subgradient algorithms. While the discontinuous subgradient of cost function is directly employed in aforementioned algorithms, which may cause vibrations of systems. The second category includes distributed smooth algorithms \cite{DCFO1}-\cite{OFW} which employed splitting method \cite{PM}. Most existing works of distributed smooth algorithms only consider one or two proximal operators \cite{CTOS3,OFW} in their algorithms. They can not directly solve the nonsmooth resource allocation optimization problem where multiple nonsmooth functions are contained in each local cost function, since summation of multiple proximable nonsmooth functions may not be proximable. \cite{Dhingra} designed a proximal augmented Lagrangian and achieved continuous-time primal-dual dynamics to solve nonsmooth optimization problems. However, more extension works are needed to solve distributed nonsmooth resource allocation optimization problems with multiple nonsmooth functions and directed graphs.

\subsection*{Contribution}
In this paper, two smooth primal-dual algorithms are proposed for a class of distributed nonsmooth convex resource allocation problems with directed graphs. A distributed estimator of the left eigenvector associated with zero eigenvalue of Laplacian matrix of the directed graph is considered in the second algorithm. The global cost function in these problems is a summation of local cost functions, and each of them consists of a smooth convex function and multiple nonsmooth convex functions. Although each  nonsmooth function is proximable, their summation might not be. Contributions of this paper are summarized as follows. 

\textbf{(i)} This paper explores a class of nonsmooth resource allocation problems with directed graphs. Compared with \cite{YP2016}-\cite{SP3}, this paper considers resource allocation problems with weight-unbalanced graphs. In contract to \cite{YWW}, smooth algorithms are designed for nonsmooth resource allocation problems with local constraints.

\textbf{(ii)} Distributed smooth primal-dual algorithms employing multi-proximal splitting are proposed in this paper. The multi-proximal splitting is used to deal with the unproximable property of the summation of nonsmooth functions and ensure smoothness of proposed algorithms.

\textbf{(iii)} A Lyapunov function and an iISS-Lyapunov function with respect to the set of equilibria are designed. Then the convergence and correctness of proposed algorithms are proved by using Lyapunov stability theory and iISS theory, which provides novel insights into analysis of the asymptotically convergent system with inputs by employing iISS theory with respect to set.

\subsection*{Organization}
The rest of this paper is organized as follows. In Section II, some basic definitions of graph theory, proximal operator and iISS theory are presented. Section III shows the nonsmooth resource allocation problem with directed graph. In Section IV, we propose two distributed multi-proximal splitting based smooth continuous-time primal-dual algorithms with and without left eigenvector estimator, respectively. Then proofs for the convergence and correctness of these algorithms are also presented. In Section V, simulations show the effectiveness of our proposed algorithm. Finally, Section VI concludes this paper.

\section{Mathematical Preliminaries}
In this section, we introduce necessary notations, definitions and preliminaries about graph theory, proximal operator and integral input-to-state stability (iISS). 

\subsection{Graph Theory}
A weighted directed graph $\mathcal{G}$ is denoted by $\mathcal{G(V,E,A)}$, where $\mathcal{V} = \lbrace 1, \dots, n \rbrace$ is a set of nodes, $\mathcal{E}$ is a set of edges, and $\mathcal{A} = [a_{ij}] \in \mathbb{R}^{n \times n}$ is a weighted adjacency matrix. An edge $e_{ij} \in \mathcal{E}$ indicates that agent $i$ can receive information from agent $j$. If $e_{ij} \in \mathcal{E}$, then $a_{ij} > 0$; otherwise, $a_{ij} = 0$. Moreover, $a_{ii} = 0, i \in \mathcal{I}$. Agent $j \in \mathcal{N}_{i}$ denotes agent $j$ is a neighbour of agent $i$. The in-degree and out-degree of agent $i$ are $d^{in}_{i} = \sum_{j=1}^{n} a_{ij}$ and $d^{out}_{i} = \sum_{j=1}^{n} a_{ji}$, respectively. The Laplacian matrix is $L_{n} = D^{in} - \mathcal{A}$, where $D^{in} \in \mathbb{R}^{n \times n}$ is diagonal with $D^{in}_{ii} = \sum^{n}_{j=1} a_{ij}$, $i \in \lbrace 1, \dots, n \rbrace$. We use $\Vert \cdot \Vert $ to indicate Euclidean norm. Let $\mathbb{R}$ denote the set of real numbers. $\mathbb{R}^{+}$ denotes the set of positive real numbers. $diag \lbrace b_{1}, \cdots, b_{n} \rbrace \in \mathbb{R}^{n \times n}$ is denoted as the diagonal matrix, whose $i$-th diagonal element is $b_{i} \in \mathbb{R}$ for $i \in \lbrace 1, \cdots, n \rbrace$. $I_{n}$ is the $n$-dimensional identity matrix. Let $\textbf{0}_{n} \in \mathbb{R}^{n}$ denote the vector of all zeros. $O_{n}$ is the $n$-dimensional null matrix, which means that every element in $O_{n}$ is zero. $(\cdot)^{T}$ denotes transpose of matrix.

\begin{lem}
(\cite{LH+HL}) Assume that graph $\mathcal{G}$ is strongly connected with the Laplacian matrix $L_{n}$. Then: 

$(1)$ There is a positive left eigenvector $h = (h_{1}, h_{2}$, $\cdots,h_{n})^{T}$ associated with the zero eigenvalue such that $h^{T}L = \textbf{0}^{T}_{n}$ and $\sum_{i = 1}^{n} h_{i} = 1$.

$(2)$ $\min_{\textbf{1}^{T}_{n}x=0} x^{T}\textbf{L} x \geq \lambda_{2}(\textbf{L}) \Vert x \Vert^{2}$, where $\textbf{L} = (HL + L^{T}H)/2$ with $H = diag(h_{1}, h_{2}, \cdots , h_{n})$ and $\lambda_{2}(\textbf{L})$ being its second smallest eigenvalue. \label{L2.1}
\end{lem}

\subsection{Proximal Operator}
Let $f(\delta)$ be a lower semi-continuous convex function for $\delta \in \mathbb{R}^{r}$. Then the proximal operator $prox_{f}[\theta]$ of $f(\delta)$ at $\theta \in \mathbb{R}^{r}$ is
\begin{equation}
prox_{f}[\theta] = \arg \min_{\delta} \lbrace f(\delta) + \frac{1}{2} \Vert \delta - \theta \Vert^{2} \rbrace.
\end{equation} 

Let $ \partial f(\delta)$ denote the subdifferential of $f(\delta)$. If $f(\delta)$ is convex, then $\partial f(\delta)$ is monotone, that is, $(\zeta_{\delta_{1}} - \zeta_{\delta_{2}})^{T}(\delta_{1} - \delta_{2}) \geq 0$ for all $\delta_{1} \in \mathbb{R}^{r}, \delta_{2} \in \mathbb{R}^{r}$, $\zeta_{\delta_{1}} \in \partial f(\delta_{1})$, and $\zeta_{\delta_{2}} \in \partial f(\delta_{2})$. $\delta =
prox_{f} [\theta]$ is equivalent to
\begin{equation}
\theta - \delta \in \partial f(\delta).
\label{Proximal Property}
\end{equation}

\subsection{Integral Input-to-State Stability with respect to set}
Consider the system
\begin{equation}\label{system}
\dot{x} = f(x,u), x(0)=x_{0}, t \geq 0,
\end{equation}
where $x \in \mathbb{R}^{n}$ and $u \in \mathbb{R}^{m}$. Inputs are measurable and locally essentially bounded functions $u: \mathbb{R}_{\geq 0} \to \mathbb{R}^{m}$, and $f: \mathbb{R}^{n}\times\mathbb{R}^{m} \to \mathbb{R}^{n}$ is assumed to be locally Lipschitz continuous. Equilibria of system \eqref{system} consist a closed set $\mathcal{M}$. For each $\xi \in \mathbb{R}^{n}$, the point-to-set distance from $\xi$ to $\mathcal{M}$ is denoted by
\begin{gather}
\Vert \xi \Vert_{\mathcal{M}} \triangleq d(\xi,\mathcal{M})=\inf\lbrace \Vert \xi -\psi \Vert, \psi \in \mathcal{M} \rbrace
\end{gather}
In particular, $\Vert \xi \Vert_{\lbrace 0 \rbrace} = \Vert \xi \Vert$. Let $\mathcal{K}$ denote the class of functions $a(x): [0,\infty) \to [0,\infty)$ which are strictly increasing, continuous and $a(0)=0$; $\mathcal{K}_{\infty}$ denotes the class of functions $a(x): [0,\infty) \to [0,\infty)$ which are a subset of $\mathcal{K}$ functions that $\lim_{x\to\infty} a(x)$ $\to \infty$; $\mathcal{L}$ is the set of functions $a(x):[0,+\infty) \to [0,+\infty)$ which are continuous, decreasing and $\lim_{x\to+\infty}a(x)=0$; $\mathcal{KL}$ is the class of functions $a(x,y):[0,\infty)^{2} \to [0,\infty)$ where $a(x,y)$ belongs to class $\mathcal{K}$ with respect to $x:[0,\infty)$ and to class $\mathcal{L}$ with respect to $y:[0,\infty)$ \cite{K}. A positive definite function $a(x):[0,\infty) \to [0,\infty)$ is one that $a(0)=0$ and $a(x)>0$ when $x>0$. A function $V(x) \in \mathbb{R}$ is semiproper if and only if for each $r$ in the range of $V(x)$, the sublevel set $\lbrace x \vert V(x) \leq r \rbrace$ is compact. A positive definite function with respect to $\mathcal{M}$ is one that is zero at $\mathcal{M}$ and positive otherwise  \cite{iISS,ISS_Lin}. A nonempty set $\mathcal{M}$ is 0-invariant for system \eqref{system} if the solution starting from $\mathcal{M}$ is defined for all $t \geq 0$ and stays in $\mathcal{M}$ when $u \equiv \textbf{0}_{m}$. System \eqref{system} is said to be forward complete if the solution $x(t,x_{0},u)$ is defined for all $t>0$ \cite{In_Lya}.

Define $DV(x) = [\frac{\partial V(x)}{\partial x}]^{T}$. Then definitions of integral input-to-state stability (iISS) and iISS-Lyapunov function with respect to a closed and 0-invariant set $\mathcal{M}$ are given below.

\begin{defn}\label{D1}
System \eqref{system} is \textbf{Integral Input-to-State Stability (iISS) with respect to a closed and 0-invariant set $\mathcal{M}$} if system \eqref{system} is forward complete and there exist functions $a_{1} \in \mathcal{K}_{\infty}$, $a_{2} \in \mathcal{KL}$ and $a_{3} \in \mathcal{K}$, such that
\begin{equation}
a_{1}(\Vert x(t,x_{0},u) \Vert_{\mathcal{M}}) \!\leq\! a_{2}(\Vert x_{0} \Vert_{\mathcal{M}}, t) + \!\!\!\int^{t}_{0}\!\!\!\! a_{3}(\Vert u(s) \Vert)ds.\label{IISS}
\end{equation}  
\end{defn}
\begin{defn} \label{D3}
A continuously differentiable function $V$ is called an \textbf{iISS-Lyapunov function with respect to a closed and 0-invariant set $\mathcal{M}$} for system \eqref{system} if system \eqref{system} is forward complete and there exist functions $a_{4}, a_{5} \in \mathcal{K}_{\infty}$ and a continuous positive definite function $a_{6}$, and $a_{7} \in \mathcal{K}$ such that 
\begin{equation}
a_{4}(\Vert x \Vert_{\mathcal{M}}) \leq V(x) \leq a_{5}(\Vert x \Vert_{\mathcal{M}}),
\end{equation}
and 
\begin{equation}
DV(x)f(x,u) \leq - a_{6}(\Vert x \Vert_{\mathcal{M}}) + a_{7}(\Vert u \Vert)
\end{equation}
for all $x \in \mathbb{R}^{n}$ and all $u \in \mathbb{R}^{m}$.
\end{defn}

Note that $V$ in Definition \ref{D3} is positive definite and proper (i.e., radially unbounded) with respect to $\mathcal{M}$. If the 0-input system $\dot{x} = f(x,\textbf{0}_{m}) \label{0system}$ is \textbf{globally asymptotically stable} (GAS) with respect to $\mathcal{M}$, the system \eqref{system} is to be said 0-GAS with respect to $\mathcal{M}$.

Similar to definitions of dissipation and zero-output dissipation in \cite{iISS2}, here we introduce concepts of dissipation and zero-output dissipation with respect to $\mathcal{M}$.

\begin{defn}\label{D2}
The system \eqref{system} with output p$: \mathbb{R}^{n} \to \mathbb{R}^{r}$ is \textbf{dissipative with respect to a closed and 0-invariant set $\mathcal{M}$} if system \eqref{system} is forward complete and there exists a continuously differentiable, proper, and positive definite function $V$ with respect to $\mathcal{M}$, together with a continuous positive definite function $a_{8}$ and a function $a_{9} \in \mathcal{K}$, such that
\begin{equation}
DV(x) f(x,u) \leq - a_{8}(\Vert p(x) \Vert) + a_{9}(\Vert u \Vert)\label{Dissipative}
\end{equation}
for all $x \in \mathbb{R}^{n}$ and all $u \in \mathbb{R}^{m}$. Moreover, if \eqref{Dissipative} holds with $p = \textbf{0}_{r}$, i.e., if there exist a proper and positive definite $V$ with respect to $\mathcal{M}$, and an $a_{9} \in \mathcal{K}$, such that
\begin{equation}
DV(x) f(x,u) \leq a_{9}(\Vert u \Vert)\label{Dissipative0}
\end{equation}
holds for all $x \in \mathbb{R}^{n}$ and all $u \in \mathbb{R}^{m}$, we say that the system \eqref{system} is \textbf{zero-output dissipative (ZOD) with respect to $\mathcal{M}$}. 
\end{defn}

Consider a system 
\begin{equation}
\dot{x}(t) = J(x(t)), x(0) = x_{0}, t \geq 0 \label{syt}
\end{equation}
where $J: \mathbb{R}^{n} \to \mathbb{R}^{n}$ is Lipschitz continuous. The following result is a special case of Theorem 3.1 in \cite{OE}.

\begin{lem} \label{Aequilibrium}
Let $\mathcal{D}$ be a compact, positive invariant set with respect to system \eqref{syt}, $V: \mathbb{R}^{n} \to \mathbb{R}$ be a continuously differentiable function, and $x(\cdot) \in \mathbb{R}^{q}$ be a solution of \eqref{syt} with $x(0) = x_{0} \in \mathcal{D}$. Assume $\dot{V}(x) \leq 0$, $\forall x \in \mathcal{D}$, and define $\mathcal{Z} = \lbrace x \in \mathcal{D}: \dot{V}(x) = 0 \rbrace$. If every point in the largest invariant subset $\mathcal{M}$ of $\bar{\mathcal{Z}}\cap\mathcal{D}$ is Lyapunov stable, where $\bar{\mathcal{Z}}$ is the closure of $\mathcal{Z} \subset \mathbb{R}^{n}$, then system \eqref{syt} converges to one of its equilibria.
\end{lem}

\section{Problem Description}
In this section, the resource allocation problem with a directed graph is formulated. We consider a network of $n$ agents with first-order dynamics, interacting over a graph $\mathcal{G}$. The nonsmooth resource allocation problem is given as
\begin{gather}
\min_{x \in \mathbb{R}^{nq}} F(x), \ \ s.t. \sum_{i=1}^{n} x_{i} = \sum_{i=1}^{n} d_{i}, \label{Problem 2}
\end{gather}
where $F(x) = \sum_{j=0}^{m} F^{j}(x) = \sum_{i=1}^{n} f_{i}(x_{i})$, $f_{i}(x_{i}) = \sum_{j=0}^{m}f^{j}_{i}(x_{i})$, and $F^{j}(x) = \sum_{i=1}^{n} f^{j}_{i}(x_{i})$, $j \in \lbrace 0, 1, \cdots$,$m \rbrace$, $m \geq 2$. Note that $x_{i} \in \mathbb{R}^{q}$ is the state of $i$-th agent and $x = [ x_{1}^{T}, x_{2}^{T}, \cdots, x_{n}^{T} ]^{T} \in \mathbb{R}^{nq}$. 

For each agent $i \in \lbrace 1, \cdots, n \rbrace$, there are $m+1$ function $f_{i}^{0}, \cdots, f_{i}^{m} : \mathbb{R}^{q} \to \mathbb{R}$, contained in the local cost function $f_{i}(x_{i}): \mathbb{R}^{q} \to \mathbb{R}$, where $f^{0}_{i}$ is a smooth convex function, $f^{j}_{i}$ is a nonsmooth convex function for $j \in \lbrace 1,\cdots,m\rbrace$. Each agent $i$ only has the information about $f_{i}^{j}$ for $j \in \lbrace 0, 1, \cdots, m \rbrace$. The constraint presented in \eqref{Problem 2} indicates that all solutions must achieve resource allocation conditions $\sum_{i=1}^{n} x_{i} = \sum_{i=1}^{n} d_{i}$. Each agent only exchanges information with its neighbours in a fully distributed manner.

Assumptions below are made for the wellposedness of the problem \eqref{Problem 2} in this section.

\begin{assum} $f^{0}_{i}$ is twice continuously differentiable and strongly convex for all $i \in \lbrace 1, \cdots, n \rbrace$, which means that there exists a constant $c > 0$ such that for agent $i$,
\begin{equation}
(\nabla f^{0}_{i}(\vartheta_{1}) - \nabla f^{0}_{i}(\vartheta_{2}))^{T}(\vartheta_{1} - \vartheta_{2}) \geq c \Vert \vartheta_{1} - \vartheta_{2} \Vert^{2},
\label{Strongly Convex}
\end{equation} 
where $\vartheta_{1} \in \mathbb{R}^{q}$, $\vartheta_{2} \in \mathbb{R}^{q}$, $\vartheta_{1} \neq \vartheta_{2}$. Without loss of generality, we assume $c > m-1$. \label{A1}
\end{assum}

\begin{assum} Each $f^{j}_{i}$ is (nonsmooth) lower semi-continuous closed proper convex functions for all $i \in \lbrace 1, \cdots, n \rbrace$, $j \in \lbrace 1, \cdots, m \rbrace$, and it is proximable. \label{A2}
\end{assum} 

\begin{assum}
The weighted graph $\mathcal{G}$ is directed and strongly connected. \label{A3}
\end{assum}

\begin{assum}
There exists at least one feasible point to problem \eqref{Problem 2}. \label{A4}
\end{assum} 

\begin{rem} The condition $c>m-1$ in Assumption \ref{A1} is mild. When $0 < c \leq m-1$, there always exists a function $f^{0'}_{i}(x) = K f^{0}_{i}(x)$ for agent $i$ with $ K >\frac{m - 1}{c}$ such that
$ (\nabla \! f^{0'}_{i}(\vartheta_{1}) - \nabla \! f^{0'}_{i}(\vartheta_{2}))^{T}\!(\vartheta_{1} - \vartheta_{2})\! \geq Kc \Vert \vartheta_{1} - \vartheta_{2} \Vert^{2} \! > \! (m - 1) \Vert \vartheta_{1} - \vartheta_{2} \Vert^{2}$. $\hfill$ $\blacklozenge$
\end{rem}

Then, we arrive at the following lemma by the Karush-Kuhn-Tucker (KKT) condition of convex optimization problems.

\begin{lem}
Under Assumptions $\ref{A1}$-$\ref{A4}$, a feasible point $x^{*} \in \mathbb{R}^{nq}$ is a solution of problem \eqref{Problem 2} if and only if there exist $x^{*} \in \mathbb{R}^{nq}$,  a constant $v_{0} \in \mathbb{R}^{q}$, and $v^{*} \in \mathbb{R}^{nq}$ such that
\begin{subequations}
\begin{align}
& \textbf{0}_{nq} \in \nabla F^{0}(x^{*}) + \sum_{j=1}^{m} \partial F^{j}(x^{*}) - v^{*}, \label{KKT1}\\
& \sum_{i=1}^{n} x_{i}^{*} = \sum_{i=1}^{n} d_{i}, v_{i}^{*} = v_{0} \text{ for } i \in \lbrace 1,\cdots,n \rbrace, \label{KKT2}
\end{align}\label{KKT}%
\end{subequations}
where $v = [ v_{1}^{T}, v_{2}^{T}, \cdots, v_{n}^{T} ]^{T}$ is the Lagrange multiplier, $\nabla F^{0}(x) = [(\nabla f^{0}_{1}(x_{1}))^{T}\!\!$, $(\nabla f^{0}_{2}(x_{2}))^{T} \!\!, \cdots$, $(\nabla f^{0}_{n}(x_{n}))^{T} ]^{T}\!\!$, and $\partial F^{j}(x) \!=\! [(\partial f^{j}_{1}(x_{1}))^{T}\!$, $(\partial f^{j}_{2}(x_{2}))^{T}, \cdots$, $(\partial f^{j}_{n}(x_{n}))^{T} ]^{T}$ for $j \in \lbrace 1, \cdots, m \rbrace$. 
\label{LKKT}
\end{lem}

The proof of Lemma \ref{LKKT} is omitted since it is a trivial extension of the proof for Theorem 3.25 in \cite{NO}.

\section{Distributed Algorithms with Multi-Proximal Operator}

The purpose of this section is to design two continuous-time distributed algorithms based on multi-proximal splitting to solve the nonsmooth resource allocation problem \eqref{Problem 2} for two cases that with known left eigenvector $h$ and with a distributed estimator of left eigenvector $h$, respectively.

In order to tackle the difficulty caused by the unproximable property of $\sum_{j=1}^{m}f^{j}_{i}(x_{i})$ for each agent $i$, here we introduce a class of auxiliary variables $z^{j}(t) \in \mathbb{R}^{nq}$ for $j \in \lbrace 1, \cdots, m-1 \rbrace$ combined with a constant parameter $\gamma \in \mathbb{R}^{+}$ such that there exist feasible points $z^{j*}$ splitting (\ref{KKT1}) as
\begin{subequations}
\begin{align}
- \nabla F^{0}(x^{*}) + v^{*} + \gamma \sum_{j=1}^{m-1}  z^{j*} \in \partial F^{m}(x^{*}), \label{z1}\\
- \gamma z^{j*} \in \partial F^{j}(x^{*}), j \in \lbrace 1,\cdots,m-1 \rbrace . \label{z2}
\end{align}
\label{z}%
\end{subequations}
According to the property \eqref{Proximal Property} of proximal operator, we can transfer \eqref{z} as 
\begin{gather}
\begin{split}
& x^{*} = Prox_{F^{m}}[x^{*} - \nabla F^{0}(x^{*}) + v^{*} + \gamma \sum_{j=1}^{m-1} z^{j*}], \\
& x^{*} = Prox_{F^{j}}[x^{*} - \gamma z^{j*}], j \in \lbrace 1,\cdots,m-1 \rbrace, \label{ABC}
\end{split}
\end{gather}
where for any $\xi = [\xi_{1}^{T}, \xi_{2}^{T}, \cdots, \xi_{n}^{T}]^{T} \in \mathbb{R}^{nq}$, $\xi_{i} \in \mathbb{R}^{q}$, $i \in \lbrace 1, \cdots, n \rbrace$, $Prox_{F^{j}}[\xi] = [(prox_{f^{j}_{1}}[\xi_{1}])^{T}, (prox_{f^{j}_{2}}[\xi_{2}])^{T}$, $\cdots, (prox_{f^{j}_{n}}[\xi_{n}])^{T}]^{T}$. $x^{*}$ and $v^{*}$ are defined in \eqref{KKT}. From \eqref{z2}, it is clear that $-\gamma z^{j*}$ is presented to estimate a subgradient in $\partial F^{j}(x^{*})$ for $j \in \lbrace 1, \cdots, m-1 \rbrace$.

\subsection{Algorithm Design with Known Left Eigenvector h}
In this subsection, we present a distributed smooth multi-proximal primal-dual algorithm for solving problem ($\ref{Problem 2}$) with the information of left eigenvector $h$. 

According to \eqref{KKT} and \eqref{z1}, we propose a smooth algorithm as
\begin{eqnarray}
\dot{x}_{i}(t) & = & prox_{f^{m}_{i}} \!\Big[\! x_{i}(t) \!-\! \nabla \!f^{0}_{i}\!(x_{i}(t)) \!+\! v_{i}(t) \!+\! \gamma\! \sum_{j=1}^{m-1}\!\! z_{i}^{j}(t) \Big] \!\!-\! x_{i}(t), \notag \\
\dot{z}_{i}^{j}(t) & = & prox_{f^{j}_{i}}[x_{i}(t) - \gamma z_{i}^{j}(t)] - x_{i}(t), \notag \\
\dot{v}_{i}(t) & = & - h_{i}^{-1}(x_{i}(t) - d_{i}) \!-\! \alpha \sum_{k \in \mathcal{N}_{i}} a_{ik} (v_{i}(t) - v_{k}(t))\! -\! w_{i}(t), \notag \\
\dot{w}_{i}(t) & = & \alpha \sum_{k \in \mathcal{N}_{i}} a_{ik} (v_{i}(t) - v_{k}(t)), \quad w_{i}(0)=\textbf{0}_{q},
\label{Algorithm 11}
\end{eqnarray}
where $t \geq 0$, $0 < \gamma < \frac{1}{m-1}$, $i \in \lbrace 1, \cdots, n \rbrace$, and $j \in \lbrace 1,\cdots,m-1 \rbrace$.

\begin{rem}
Because all proximal operators $prox_{f^{j}_{i}}(\cdot)$ for $i \in \lbrace 1, \cdots, n \rbrace$ and $ j \in \lbrace 1, \cdots, m-1 \rbrace$ are continuous and nonexpansive, the proposed algorithm \eqref{Algorithm 11} is locally Lipschitz continuous even though each $f^{j}_{i}(x_{i})$ in problem \eqref{Problem 2} is nonsmooth, which means that the \textbf{smoothness} of algorithm \eqref{Algorithm 11} is guaranteed. $\hfill$ $\blacklozenge$
\end{rem} 

Algorithm \eqref{Algorithm 11} can be written in a compact form as
\begin{subequations}
\begin{align}
\dot{x}(t) = & Prox_{F^{m}} \Big[ x(t) \!\!-\!\! \nabla F^{0}(x(t)) \!\!+\!\! v(t) \notag \\
& +\!\! \gamma \sum_{j = 1}^{m-1} z^{j}(t) \Big] \!\!-\!\! x(t), \label{another proximal} \\
\dot{z}^{j}(t) = & Prox_{F^{j}}[x(t) - \gamma z^{j}(t)] - x(t),  \label{estimator} \\
\dot{v}(t) = & - H^{-1}_{nq}(x(t) - d) - \alpha L_{nq} v(t) - w(t), \\
\dot{w}(t) = & \alpha L_{nq} v(t), \quad w(0)=\textbf{0}_{nq}, \label{lagrangian}
\end{align}\label{Algorithm 21}%
\end{subequations}
where $j \in \lbrace 1,\cdots,m-1 \rbrace$, $H_{nq} = diag \lbrace h_{1}, \cdots, h_{n} \rbrace \otimes I_{q}$, $d=[d_{1}^{T}, \cdots, \d_{n}^{T}]^{T} \in \mathbb{R}^{nq}$, and $L_{nq}= L_{n} \otimes I_{q}$. The matrix $L_{n} \otimes I_{q}$ is the Kronecker product of matrices $L_{n}$ and $I_{q}$.

\begin{rem}
From \eqref{estimator}, it is shown that $-\gamma z^{j}$ is the proximal-based estimator of a subgradient in $\partial F^{j}(x)$ for $j \in  \lbrace 1, \cdots, m-1 \rbrace$. With the help of estimator $- \gamma z^{j}$, the corresponding proximal operator \eqref{another proximal}, which employs the information of $-\gamma z^{j}$ instead of $\partial F^{j}(x)$ for $j \in  \lbrace 1, \cdots, m-1 \rbrace$, is presented to tackle the difficulty caused by the unproximable property of $\sum_{j=1}^{m-1} F^{j}(x)$. The scheme combined by \eqref{another proximal} and \eqref{estimator} is called the \textbf{multi-proximal splitting}, which may be viewed as an extension of three operator splitting shown in \cite{OFW}. $\hfill$ $\blacklozenge$
\end{rem} 

\begin{lem}
Under Assumptions $\ref{A1}$-$\ref{A4}$, if $(x^{*}$, $z^{*}$, $v^{*}$, $w^{*}) \in (\mathbb{R}^{nq}, \mathbb{R}^{(m-1)nq}, \mathbb{R}^{nq}, \mathbb{R}^{nq})$ is an equilibrium of algorithm \eqref{Algorithm 21} and $(\textbf{1}_{n}\otimes I_{q})^{T}H_{nq}w^{*}=\textbf{0}_{q}$, then $x^{*}$ is a solution of problem \eqref{Problem 2}, where $z = [(z^{1})^{T}, \cdots, (z^{2})^{T}]^{T}$. \label{KA1}
\end{lem}

\begin{pf}
If $(x^{*}, z^{*}, v^{*}, w^{*})$ is an equilibrium of algorithm \eqref{Algorithm 21}, then according to the property \eqref{Proximal Property} of proximal operator and algorithm \eqref{Algorithm 21}, it yields that for $j \in \lbrace 1,\cdots,m -1 \rbrace$, 
\begin{subequations}
\begin{align}
- \nabla F^{0}(x^{*}) + v^{*} + \gamma \sum_{j=1}^{m-1} z^{j*} \!\!\in & \partial F^{m}(x^{*}), 
\label{P11} \\
- \gamma z^{j*} \!\!\in & \partial F^{j}(x^{*}), \label{P21} \\
- H^{-1}_{nq}(x^{*} - d) - \alpha L_{nq}v^{*} - w^{*} \!\!= & \textbf{0}_{nq}, 
\label{P31} \\
\alpha L_{nq} v^{*} \!\!= & \textbf{0}_{nq}. \label{P41} 
\end{align}
\end{subequations}

From \eqref{P11}, \eqref{P21} and \eqref{P41}, there exists a $v^{0} \in \mathbb{R}^{q}$ such that
\begin{gather}\label{L4A}
\begin{split}
\textbf{0}_{nq} \in & - \nabla F^{0}(x^{*}) - \sum_{j=1}^{m-1} \partial F^{j}(x^{*}) + v^{*}, \\
v^{*} = & \textbf{1}_{n} \otimes v^{0}.
\end{split}
\end{gather}


Summing \eqref{P31} and \eqref{P41} yields that $-(x^{*} - d) - H_{nq}w^{*} = \textbf{0}_{nq}$, which means that
\begin{equation}\label{L4B}
\sum_{i = 1}^{n} (x^{*}_{i} \!\!-\! d_{i})\! =\! -\!\! \sum_{i = 1}^{n} \! h_{i}I_{q}w^{*}_{i}\!\! =\!\! -(\textbf{1}_{n} \otimes \! I_{q})^{T}\!H_{nq}w^{*} \!\! = \textbf{0}_{q}.
\end{equation}

Considering \eqref{L4A} together with \eqref{L4B} and according to Lemma $\ref{LKKT}$, $x^{*}$ is a solution of problem \eqref{Problem 2}. $\hfill$ $\blacksquare$
\end{pf}

Then we state the convergence result of the proposed distributed algorithm \eqref{Algorithm 21}. Let $(x^{*}, z^{*}, v^{*}, w^{*})$ be an equilibrium of algorithm \eqref{Algorithm 21}. Define a Lyapunov candidate $V(x,z,v,w) = V_{1}(x,z) + V_{2}(x) + V_{3}(v,w)$, where
\begin{eqnarray}\label{VV123}
& & V_{1}(x,z) = (\eta \! + \! 1) [\frac{1}{2} \Vert \bar{x}^{*} \Vert^{2} \!+ \! \frac{1}{2} \gamma \sum_{j=1}^{m-1} (\Vert \bar{z}^{j*} \Vert^{2}  \! - \! 2 (\bar{x}^{*})^{T}\bar{z}^{j*})], \notag \\
& & V_{2}(x) = (\eta \! + \! 1) [F^{0}(x) \! -  \! F^{0}(x^{*}) - (\bar{x}^{*})^{T} \nabla F^{0}(x^{*})],  \\
& & V_{3}(v,w) = \frac{\eta}{2} (\bar{v}^{*})^{T}H_{nq}\bar{v}^{*} + \frac{1}{2} (\bar{v}^{*} + \bar{w}^{*})^{T}H_{nq}(\bar{v}^{*} + \bar{w}^{*}), \notag 
\end{eqnarray}
and $\eta > 0$, $\bar{x}^{*} \triangleq x - x^{*}$, $\bar{z}^{j*} \triangleq z^{j} - z^{j*}$, $\bar{v}^{*} \triangleq v - v^{*}$, $\bar{w}^{*} \triangleq w - w^{*}$.

By analysing the convergence of \eqref{Algorithm 21}, the main theorem of this subsection is obtained as below.

\begin{thm}
Consider algorithm \eqref{Algorithm 21}. Suppose Assumptions \ref{A1}-\ref{A4} hold. If following inequalities
\begin{gather}
\alpha  >  \frac{(\eta + 1)^{2}}{\eta \lambda_{2}(\textbf{L}_{nq})}, \text{ } \eta > \max \lbrace \frac{1}{b_{2}h^{*}} - 1, 0 \rbrace \label{IE}
\end{gather}
hold, where $b_{2} = c - \frac{1}{2}(1 + \gamma)(m - 1) \frac{1}{\beta}$, $\frac{(1 + \gamma)(m - 1)}{2c} < \beta < \frac{2}{1 + \gamma}$, $h^{*} = \min_{i \in \mathcal{I}} \lbrace h_{1}, \cdots, h_{n} \rbrace$, then the trajectory of $x(t)$ converges, and $\lim_{t \to \infty} x(t)$ is the solution of problem \eqref{Problem 2}. \label{Thm1}
\end{thm}
\begin{pf}
It can be easily verified that $V(x^{*}, z^{*}, v^{*}, w^{*})$ $= 0$. Next, we will show that $V(x,z,v,w) > 0$ for all $(x,z,v,w) \neq (x^{*}, z^{*}, v^{*}, w^{*})$.

Since all $f^{0}_{i}(x)$ for $i \in \lbrace 1,\cdots, n \rbrace$ are convex, then $F^{0}(x) - F^{0}(x^{*}) - (\bar{x}^{*})^{T} \nabla F^{0}(x^{*}) \geq 0$. Hence $ V_{2}(x)  \geq  0 $.

Since $0 < \gamma < \frac{1}{m-1}$, 
\begin{eqnarray}
V_{1}(x,z) & = & \frac{\eta+1}{2} \sum_{j=1}^{m-1} \left[ \Vert (\frac{1}{m-1})^{\frac{1}{2}}\bar{x}^{*} - \gamma (m-1)^{\frac{1}{2}} \bar{z}^{j*}\Vert^{2}\right. \notag \\
& & + \gamma (1 - \gamma(m-1))\Vert \bar{z}^{j*}\Vert^{2} \biggr] \geq 0. \label{V1}
\end{eqnarray}

Since $V_{2}(x) \geq 0$, $V(x,z,v,w) \geq V_{1}(x,z) + V_{3}(v,w) \geq 0 $. Clearly $V(x,z,v,w)$ is positive definite, radically unbounded, $V(x,z,v,w) \geq 0$ and is zero if and only if $(x,z,v,w) = (x^{*},z^{*},v^{*},w^{*})$. 

It follows from algorithm \eqref{Algorithm 21} that
\begin{gather}\label{x-x*}
\begin{split}
x + \dot{x} = & Prox_{F^{m}}[x - \nabla F^{0}(x) + v + \sum_{j=1}^{m-1} \gamma z^{j}], \\
x^{*} = & Prox_{F^{m}} [x^{*} \! - \! \nabla F^{0}(x^{*}) + v^{*} \! + \! \sum_{j=1}^{m-1} \gamma z^{j*}], \\
x + \dot{z}^{j} = & Prox_{F^{j}}[x - \gamma z^{j}], j \in \lbrace 1, \cdots, m-1 \rbrace, \\
x^{*} = & Prox_{F^{j}} [x^{*} - \gamma z^{j*}], j \in \lbrace 1, \cdots, m-1 \rbrace. 
\end{split}
\end{gather}

Since $f^{j}_{i}(\cdot)$ is convex, $\partial f^{j}_{i}(\cdot)$ is monotone for agent $i \in \lbrace 1, \cdots, n \rbrace$, where $j \in \lbrace 1, \cdots, m-1 \rbrace$. According to the property \eqref{Proximal Property} of proximal operator, it follows from \eqref{x-x*} that for $j \in \lbrace 1, \cdots, m-1 \rbrace$,
\begin{gather}
\begin{split}
(\gamma \sum_{j=1}^{m-1} \bar{z}^{j*}\!\! - \!\!\nabla F^{0}(\bar{x}^{*}) + \bar{v}^{*} \!\!-\! \dot{x})^{T} (\bar{x}^{*}\!\! + \!\dot{x}) \geq & 0,\\
(-\gamma \bar{z}^{j*} \!\!-\!\! \dot{z}^{j})^{T}(\bar{x}^{*} \!\!+\! \dot{z}^{j}) \geq & 0,
\end{split}\label{Prox_pro}
\end{gather}
where $\nabla F^{0}(\bar{x}^{*}) \triangleq \nabla F^{0}(x) - \nabla F^{0}(x^{*})$.

From \eqref{Prox_pro}, it can be shown that for $j \in \lbrace 1, \cdots, m-1 \rbrace$, 
\begin{gather}\label{Proxmal relation1}
\begin{split}
& \gamma \sum_{j=1}^{m-1} [(\bar{z}^{j*})^{T}\! \bar{x}^{*}] \!\!-\!\! (\nabla \! F^{0}\! (\bar{x}^{*}))^{\!T}\bar{x}^{*} \!\!+\!\! (\bar{v}^{*})^{\!T} \bar{x}^{*} \!+\! (\bar{v}^{*}\!)^{T} \! \dot{x} \\
+ & \gamma \! \sum_{j=1}^{m-1} \! [(\bar{z}^{j*})^{T} \dot{x}] \!-\! (\nabla \! F^{0}\!(\bar{x}^{*}))^{T} \!\! \dot{x} \!-\! (\bar{x}^{*})^{T} \!\!\dot{x} \!-\! \Vert \dot{x} \Vert^{2} \! \geq \! 0, 
\end{split}
\end{gather}
and
\begin{gather}\label{Proxmal relation2}
- \gamma (\bar{z}^{j*})^{T} \bar{x}^{*} - (\bar{x}^{*})^{T} \dot{z}_{j} - \gamma (\bar{z}^{j*})^{T} \dot{z}_{j} - \Vert \dot{z}_{j} \Vert^{2} \geq 0.
\end{gather}

The derivative of Lyapunov candidate $V(x,z,v,w)$ along the trajectory of algorithm \eqref{Algorithm 21} satisfies
\begin{gather}
\begin{split}
& \dot{V}(x,z,v,w) \\
= & (\eta+1)(\bar{x}^{*})^{T} \! \dot{x} \! + \! \gamma (\eta+1)\sum_{j=1}^{m-1} (\bar{z}^{j*})^{T} \! \dot{z} \\
& \! - \! \gamma (\eta+1)\sum_{j=1}^{m-1} ((\bar{x}^{*})^{T} \! \dot{z}^{j} \! + \! (\bar{z}^{j*})^{T} \! \dot{x}) \\
& + (\eta+1)(\nabla F^{0}(\bar{x}^{*}))^{T} \dot{x} + \dot{V}_{3}(v,w), \label{dot_V1}
\end{split}
\end{gather}
where
\begin{gather}\label{DV3}
\begin{split}
& \dot{V}_{3}(v,w)\\
\leq & - (\eta + 1)(\bar{v}^{*})^{T}\bar{x}^{*} - (\bar{w}^{*})^{T}\bar{x}^{*} - (\bar{w}^{*})^{T}H_{nq}\bar{w}^{*} \\
& - (\eta + 1)(\bar{v}^{*})^{T}H_{nq}\bar{w}^{*} - \alpha \eta (\bar{v}^{*})^{T}\textbf{L}_{nq}\bar{v}^{*}, 
\end{split}
\end{gather}
and $\textbf{L}_{nq} = (H_{nq}L_{nq} + L_{nq}^{T}H_{nq})/2$.

According to \eqref{Proxmal relation1}-\eqref{DV3}, it follows that 
\begin{gather}
\begin{split}
& \dot{V}(x,z,v,w) \\
\leq & - (\eta + 1) \Vert \dot{x} \Vert^{2} - (\eta + 1) \sum_{j=1}^{m-1} \Vert \dot{z}^{j} \Vert^{2}  \\
& - (\eta + 1)(1 + \gamma) \sum_{j=1}^{m-1}(\bar{x}^{*})^{T} \dot{z}^{j} - (\bar{w}^{*})^{T}H_{nq}\bar{w}^{*}  \\
& - (\eta + 1) (\nabla F^{0}(\bar{x}^{*}))^{T} \! \bar{x}^{*} + (\eta + 1)(\bar{v}^{*})^{T}\dot{x}  \\
&  - (\bar{w}^{*})^{T}\bar{x}^{*} - \! \alpha \eta (\bar{v}^{*})^{T}\textbf{L}_{nq}\bar{v}^{*} - (\eta + 1)(\bar{v}^{*})^{T}\bar{w}^{*}.
\end{split}
\end{gather}

Then according to Assumption \ref{A1}, there exists a parameter $\beta > 0$ such that 
\begin{gather}
\begin{split}
& (1 + \gamma) \sum_{j=1}^{m-1} (\bar{x}^{*})^{T}\!\dot{z}^{j} \\
\geq & -\! \frac{1}{2}(1\! + \gamma) \beta \sum_{j=1}^{m-1} \Vert \dot{z}^{j} \Vert^{2} - \frac{(1 + \gamma)(m-1)}{2\beta} \Vert \bar{x}^{*} \Vert^{2}.
\end{split}
\end{gather}

Hence we have the conclusion that
\begin{gather}\label{DV}
\begin{split}
& \dot{V}(x,z,v,w) \\
\leq & -\! (\eta \! + \! 1) \Vert \dot{x} \Vert^{2} - (\eta + 1) b_{1} \sum_{j=1}^{m-1} \Vert \dot{z}^{j} \Vert^{2} - (\bar{w}^{*})^{T}\bar{x}^{*}  \\ 
& -\! \alpha \eta (\bar{v}^{*}\!)^{T}\textbf{L}_{nq}\bar{v}^{*} \!\!-\! (\eta \! + \! 1)b_{2}\Vert \bar{x}^{*} \Vert^{2} \!\!-\! (\bar{w}^{*}\!)^{T} \! H_{nq}\bar{w}^{*} \\
& + (\eta + 1)(\bar{v}^{*})^{T}\dot{x} - (\eta + 1)(\bar{v}^{*})^{T}H_{nq}\bar{w}^{*},
\end{split}
\end{gather}
where $b_{1} = 1 - \frac{1}{2}(1 + \gamma) \beta$ and $b_{2} = c - \frac{1}{2}(1 + \gamma)(m-1) \frac{1}{\beta}$.

In order to illustrate that there always exists a $\beta > 0$ such that $b_{1} > 0$ and $b_{2} > 0$, here we define a function $B(\gamma)$ of $\gamma$ that $B(\gamma) = \frac{2}{\gamma + 1} - \frac{(\gamma + 1)(m - 1)}{2c}$. The derivative of $B(\gamma)$ is shown as
\begin{gather}
\frac{d B(\gamma)}{d \gamma} = -\frac{2}{(\gamma + 1)^{2}} - \frac{m - 1}{2c} < 0. \label{Beta}
\end{gather}

Note that $0 < \gamma < \frac{1}{m - 1} \leq 1$ and $c > m - 1$. According to  \eqref{Beta}, we have $B_{min}(\gamma) > B(1) = 1 - \frac{m - 1}{c} > 0$. 

As the result, there exists a $\beta$ such that 
\begin{equation}
\frac{(1 + \gamma)(m - 1)}{2c} < \beta < \frac{2}{1 + \gamma}, \label{Beta1}
\end{equation}
which means that $b_{1} = 1 - \frac{1}{2}(1 + \gamma) \beta > 0$, and $b_{2} = c - \frac{1}{2}(1 + \gamma)(m - 1) \frac{1}{\beta} > 0$.

In light of the above analysis and using the
inequality $x^{T}y \leq \frac{1}{2\tau} \Vert x \Vert^{2} + \frac{\tau}{2} \Vert y \Vert^{2}$, equation \eqref{DV} can be written as
\begin{gather}
\begin{split}
\dot{V}(x,z,v,w) \leq & \! - \epsilon_{1} \Vert \dot{x} \Vert^{2} \! - \epsilon_{2} \sum_{j=1}^{m-1} \Vert \dot{z}^{j} \Vert^{2} - \epsilon_{3}\Vert \bar{x}^{*} \Vert^{2} \\
& \! - \epsilon_{4} \Vert \bar{v}^{*} \Vert^{2} \! - \epsilon_{5} (\bar{w}^{*})^{T}H_{nq}\bar{w}^{*},
\end{split}
\end{gather}
where $\epsilon_{1} = \eta + \frac{1}{2}$, $\epsilon_{2} = (\eta + 1) b_{1}$, $\epsilon_{3} = (\eta + 1)b_{2}- \frac{1}{h^{*}}$, $\epsilon_{4} = \alpha \eta \lambda_{2}(\textbf{L}_{nq}) - (\eta + 1)^{2}$, and $\epsilon_{5} = \frac{1}{4}$.

According to \eqref{IE}, it follows that $\epsilon_{k} > 0$ for $k \in \lbrace 1,2,3,4,5 \rbrace$. Additionally, since $V(x,z,v,w)$ is positive-definite, radically unbounded, lower bounded, $(x^{*},z^{*},$ $v^{*}, w^{*})$ is Lyapunov stable. It follows from the LaSalle invariant principle and Lemma \ref{Aequilibrium} that $(x(t), z(t), v(t), w(t))$ converges to an equilibrium of algorithm \eqref{Algorithm 21} in the largest invariant set $\mathcal{M}$ in $ E = \lbrace (x,z,v,w) \vert x=x^{*}, v=v^{*},w=w^{*}, -\gamma z^{j} \in \partial F^{j}(x^{*})$ for $j \in \lbrace 1, \cdots, m-1\rbrace \rbrace$. Since $(h \otimes I_{q})^{T} L_{nq} = (\textbf{0}_{n} \otimes I_{q})^{T}$, $\sum_{i = 1}^{n} h_{i}I_{q} \dot{w}_{i}(t) = (\textbf{1}_{n}\otimes I_{q})^{T}H_{nq} L_{nq} v(t) = \textbf{0}_{q} $. With $w(0) = \textbf{0}_{nq}$, it shows that $\sum_{i = 1}^{n} h_{i}I_{q} w_{i}(t) = (\textbf{1}_{n}\otimes I_{q})^{T}H_{nq}w^{*} = \textbf{0}_{q}$. According to Lemma \ref{KA1}, $x^{*}$ is a solution of problem \eqref{Problem 2}. $\hfill$ $\blacksquare$ 
\end{pf}



\subsection{Algorithm Design with Distributed Estimator of Left Eigenvector h}
However, the left eigenvector $h$ corresponding to $\lambda_{1}(L_{nq})=0$ may not be known by any single agent, since $h$ is a global variable for multi-agent systems. In this subsection, we present a distributed smooth multi-proximal primal-dual algorithm for solving the problem \eqref{Problem 2} with a distributed estimator of left eigenvector $h$. 

Similar to algorithm \eqref{Algorithm 21}, according to \eqref{KKT} and \eqref{z1}, we propose a smooth algorithm as
\begin{eqnarray}\label{Algorithm 2}
\dot{x}(t) & = & Prox_{F^{m}}[x(t) \!\!-\!\! \nabla F^{0}(x(t)) \!\!+\!\! v(t) \!\!+\!\! \gamma \sum_{j = 1}^{m-1} z^{j}(t)] \!\!-\!\! x(t), \notag \\
\dot{z}^{j}(t) & = & Prox_{F^{j}}[x(t) - \gamma z^{j}(t)] - x(t), \notag \\
\dot{v}(t) & = & - Y^{-1}(t)(x(t) - d) - \alpha L_{nq} v(t) - w(t), \\
\dot{w}(t) & = & \alpha L_{nq} v(t), \quad w(0)=\textbf{0}_{nq}, \notag \\
\dot{y}(t) & = & - L_{nn} y(t), \quad y(0)=[I_{n}^{1},\cdots,I_{n}^{n}]^{T} \in \mathbb{R}^{nn}, \notag
\end{eqnarray}
where $j \in \lbrace 1, \cdots, m-1 \rbrace$, $Y = diag \lbrace y_{1}^{1}, \cdots, y_{n}^{n} \rbrace \otimes I_{q}$, $L_{nn}=L_{n}\otimes I_{n}$, and $I_{n}^{i}$ is the $i$-th row of $I_{n}$.

\begin{rem}
When the directed graph $\mathcal{G}$ of problem \eqref{Problem 2} is weight-balanced, it follows that $h_{i}=h_{j}, i,j \in \mathcal{V}$. While $\mathcal{G}$ is usually weight-unbalanced, hence a distributed estimator of $h$ is required for problem \eqref{Problem 2}. Variable $y$ in algorithm \eqref{Algorithm 2} is designed to obtain the estimated value of $h$. Lemma \ref{KA2} combined with Theorem 2 will show that $y_{i}^{i*} = h_{i}$, where $y_{i}^{i*} = \lim\limits_{t \to \infty} y_{i}^{i}(t)$ for $i \in \lbrace 1, \cdots, n \rbrace$.
\end{rem}

\begin{lem}
Under Assumptions $\ref{A1}$-$\ref{A4}$, if $(x^{*}, z^{*}, v^{*},$ $w^{*}, y^{*}) \in (\mathbb{R}^{nq}, \mathbb{R}^{(m-1)nq}, \mathbb{R}^{nq}, \mathbb{R}^{nq}, \mathbb{R}^{nn})$ is an equilibrium of algorithm \eqref{Algorithm 2}, $(\textbf{1}_{n}\otimes I_{q})^{T}H_{nq}w^{*}=\textbf{0}_{q}$ and $y^{*} = \textbf{1}_{n} \otimes h$, then $x^{*}$ is a solution of problem \eqref{Problem 2}.  \label{KA2}
\end{lem}

\begin{pf}
If $(x^{*}, z^{*}, v^{*}, w^{*}, y^{*})$ is an equilibrium of algorithm \eqref{Algorithm 2}, similar to the proof of Lemma \ref{KA1}, it can be shown that there exists a $v^{0} \in \mathbb{R}^{q}$ such that
\begin{gather}
\begin{split}
\textbf{0}_{nq} \in & - \nabla F^{0}(x^{*}) - \sum_{j=1}^{m-1} \partial F^{j}(x^{*}) + v^{*}, \\
v^{*} = & \textbf{1}_{n} \otimes v^{0},
\end{split}
\end{gather}
and
\begin{subequations}\label{PPY}
\begin{align}
- Y^{-1*}(x^{*} - d) - \alpha L_{nq}v^{*} - w^{*} = & \textbf{0}_{nq},
\label{P311} \\
\alpha L_{nq} v^{*} = \textbf{0}_{nq}, L_{nn} y^{*} = & \textbf{0}_{nq}. \label{P411} 
\end{align}
\end{subequations}

%

Adding \eqref{P311} and \eqref{P411} yields that $-(x^{*} - d) - Y^{*}w^{*} = \textbf{0}_{nq}$, which means that
\begin{gather}
\begin{split}
 \sum_{i = 1}^{n} (x^{*}_{i} - d_{i}) = & -\sum_{i = 1}^{n} y_{i}^{i*}I_{q} w_{i}(t) = -\sum_{i = 1}^{n} h_{i}I_{q}w^{*}_{i} \\
= & -(\textbf{1}_{n}\otimes I_{q})^{T}H_{nq}w^{*} = \textbf{0}_{q}, 
\end{split}
\end{gather}
where $y_{i}^{i*}$ is the $[(i-1)q+i]$-th element of $y^{*}$. According to Lemma $\ref{LKKT}$, $x^{*}$ is a solution of problem \eqref{Problem 2}. $\hfill$ $\blacksquare$
\end{pf}

Next, we will state the convergence result of the proposed distributed algorithm \eqref{Algorithm 2}. 

Firstly, some lemmas should be given to obtain the final result.

\begin{lem}\label{Diiss}
Assume system \eqref{system} can be written as
\begin{equation}\label{system+}
\dot{x} = f(x,u) = g(x) + u.  
\end{equation}
If system \eqref{system+} is forward complete, 0-GAS with respect to a closed and 0-invariant set $\mathcal{M}$, ZOD with respect to $\mathcal{M}$ with a positive definite function $W_{1}$ that 
\begin{gather}\label{Dissipative0+}
\begin{split}
a_{10}(\Vert x \Vert_{\mathcal{M}}) \leq W_{1}(x) \leq a_{11}(\Vert x \Vert_{\mathcal{M}})\\
DW_{1}(x) f(x,u) \leq a_{12}(\Vert u(t) \Vert),
\end{split}
\end{gather}
for $a_{10}, a_{11} \in \mathcal{K}_{\infty}$ and $a_{12} \in \mathcal{K}$, then system \eqref{system+} is iISS with respect to $\mathcal{M}$ with $a_{13} \in \mathcal{KL}$ such that
\begin{eqnarray}
& & a_{10}(\Vert x(t,x_{0},u) \Vert_{\mathcal{M}}) \notag \\
& \leq & a_{13}(\Vert x_{0} \Vert_{\mathcal{M}}, t) + \int^{t}_{0} 2 \left( a_{12}(\Vert u(s) \Vert) + \Vert u(s) \Vert \right) ds. \label{IISS+}
\end{eqnarray} 
Moreover, if $a_{12}(\Vert u(t) \Vert) = k \Vert u(t) \Vert^{2}$, where $k \in \mathbb{R}^{+}$, $u(t)$ is exponentially convergent to zero, then system \eqref{system+} converges to $\mathcal{M}$. 
\end{lem}
\begin{pf}
If system \eqref{system+} is forward complete and 0-GAS with respect to $\mathcal{M}$, then by Theorem 2.8 and Remark 4.1 in \cite{In_Lya}, there exists a smooth function $W_{2}: \mathbb{R}^{n} \to \mathbb{R}$ and functions $a_{14},a_{15},a_{16} \in \mathcal{K}_{\infty}$ such that
\begin{gather} \label{CL+}
\begin{split}
a_{14}(\Vert x \Vert_{\mathcal{M}}) \leq W_{2}(x) \leq a_{15}(\Vert x \Vert_{\mathcal{M}})\\
DW_{2}(x) f(x,0) \leq - a_{16}(\Vert x \Vert_{\mathcal{M}}).
\end{split}
\end{gather}
Then according to \eqref{CL+}, proof of Lemma IV.10 and Proposition II.5 in \cite{iISS2}, there exists an iISS Lyapunov function $W_{3}$ with respect to $\mathcal{M}$ such that
\begin{equation}
W_{3}(x) = W_{1}(x) + \pi (W_{2}(x)), \label{CLF++}
\end{equation}
where $\pi(r) \triangleq \int^{r}_{0} \frac{ds}{1+\kappa(a_{14}^{-1}(s))}$, and $\kappa(r) \triangleq r$ $+ \max_{ \Vert x \Vert_{\mathcal{M}} \leq r}$ $ \lbrace \Vert DW_{2}(x) \Vert\rbrace$.

From \eqref{system+}, \eqref{Dissipative0+} and \eqref{CLF++}, we have the conclusion that
\begin{gather}\label{Dissipative0++}
\begin{split}
& DW_{3}(x) f(x,u) \\
\leq & - \rho(W_{3}(x)) + (\Vert u(t) \Vert + a_{12}(\Vert u(t) \Vert)),
\end{split}
\end{gather}
where $\rho$ is a positive definite function. Then according to \eqref{Dissipative0++} and Corollary IV.3 in \cite{iISS2}, there exist an $a_{17} \in \mathcal{KL}$ such that 
\begin{gather}
\begin{split}
& W_{3}(x(t)) \\
\leq & a_{17}(W_{3}(x_{0}),t) + \int^{t}_{0} \! \! 2 (\Vert u(\tau) \Vert + a_{12}(\Vert u(\tau) \Vert)) d\tau.
\end{split}
\end{gather}
Since $a_{10}(\Vert x \Vert_{\mathcal{M}}) \leq W_{1}(x) \leq W_{3}(x) \leq W_{1}(x)+W_{2}(x)$, equation \eqref{IISS+} holds. 

Let $U(t)=\int^{\infty}_{t} 2(k\Vert u(\tau) \Vert^{2} + \Vert u(\tau) \Vert )d\tau$ for $t \geq 0$. Since $u(t) $ is exponentially convergent to zero, $U(t) \leq M_{U}$ for a $M_{U} \in \mathbb{R}^{+}$, $U(t)$ is decreasing, and $\lim_{t\to\infty}U(t)=0$. From \eqref{IISS+}, for $t \geq 0$, it follows that
\begin{gather}
\begin{split}
\Vert x(t) \Vert_{\mathcal{M}} \leq a_{10}^{-1}(a_{13}(\Vert x(0) \Vert_{\mathcal{M}},0)+M_{U}) \triangleq M_{X}.
\end{split}
\end{gather}
For any $\varepsilon > 0$, choose $T_{U} \geq 0$ and $T_{X} \geq 0$ such that $U(T_{U}) \leq a_{10}(\varepsilon)/2$ and $a_{13}(M_{X},T_{X}) \leq a_{10}(\varepsilon)/2$. Let $T \triangleq T_{X} + T_{U}$. Then from \eqref{IISS+}, for any $t \geq T$,
\begin{gather}
\begin{split}
& a_{10}(\Vert x(t) \Vert_{\mathcal{M}})\\
\leq & a_{13}(\Vert x(T_{U}) \Vert_{\!\mathcal{M}},t\!-\!T_{U})\!+\!\!\!\int^{t}_{T_{U}}\!\! \!\!2(k\Vert u(\tau) \Vert^{2}\!\! + \!\!\Vert u(\tau) \Vert )d\tau \\
\leq & a_{13}(M_{X},T_{X}+(t-T))+U(T_{U})\\
\leq & a_{13}(M_{X},T_{X})+U(T_{U}) \leq a_{10}(\varepsilon),
\end{split}
\end{gather}
which means that $\Vert x(t) \Vert_{\mathcal{M}} \leq \varepsilon$ for all $t \geq T$. Hence system \eqref{system+} converges to $\mathcal{M}$. $\hfill$ $\blacksquare$

%
%
\end{pf}

Let $\mathcal{M}_{Yj} = \lbrace [(\varphi_{1}x^{*} - \varphi_{2}z^{j*})^{T}, (\varphi_{3}z^{j*})^{T}]^{T} \vert (x^{*},z^{*},v^{*},w^{*},$ $y^{*}) $ $\in \mathcal{M}_{Y}\rbrace$ for $j \in \lbrace 1, \cdots, m-1 \rbrace$ and $\varphi_{k} \in \mathbb{R}$ for $k \in \lbrace 1, 2, 3 \rbrace$, where $\mathcal{M}_{Y}$ is the largest invariant set in $ E = \lbrace (x,z,v,w,y) \vert x=x^{*}, v=v^{*},w=w^{*}, y=y^{*}, -\gamma z^{j} \in \partial F^{j}(x^{*})$ for $j \in \lbrace 1, \cdots, m-1\rbrace \rbrace$. Then it follows that for any $\xi \in \mathbb{R}^{2nq}$, $\xi \in \mathcal{M}_{Yj}$ if and only if $\xi = [(\varphi_{1}x^{*} - \varphi_{2}z^{j*})^{T}, (\varphi_{3}z^{j*})^{T}]^{T}$ for $(x^{*},z^{*},v^{*},w^{*},y^{*}) \in \mathcal{M}_{Y}$ and $j \in \lbrace 1, \cdots, m-1 \rbrace$.

\begin{lem}\label{x-z}
Consider algorithm \eqref{Algorithm 2}. For $\xi$ and $\mathcal{M}_{Yj}$ with $j \in \lbrace 1, \cdots, m-1 \rbrace$, it follows that
\begin{itemize}
  \item [$(i)$] 
For each $\xi \in \mathbb{R}^{2nq}$, there exists a unique $x \in \mathbb{R}^{nq}$ and $z^{j} \in \mathbb{R}^{nq}$ such that $\xi = [(\varphi_{1}x - \varphi_{2}z^{j})^{T}, (\varphi_{3}z^{j})^{T}]^{T}$.
  \item [$(ii)$]
Let $P_{\mathcal{M}_{Yj}}(\xi) \triangleq \arg \min_{\psi \in \mathcal{M}_{Yj}} \lbrace \Vert \xi - \psi \Vert^{2} \rbrace$. Then $P_{\mathcal{M}_{Yj}}(\xi)=[(\varphi_{1}x^{*} - \varphi_{2}z^{j*})^{T}, (\varphi_{3}z^{j*})^{T}]^{T}$ for some $(x^{*},z^{*},v^{*},w^{*},y^{*}) \in \mathcal{M}_{Y}$.
  \item [$(iii)$]
For each $\xi \in \mathbb{R}^{2nq}$, there exists an $(x^{*},z^{*},v^{*},w^{*},y^{*}) \in \mathcal{M}_{Y}$ such that $\Vert \xi \Vert^{2}_{\mathcal{M}_{Yj}} = \Vert \varphi_{1} \bar{x}^{*} - \varphi_{2} \bar{z}^{j*} \Vert^{2} + \Vert \varphi_{3} \bar{z}^{j*} \Vert^{2}$. 
  \item [$(iv)$]
Let $V(\xi) = \frac{1}{2}\Vert \xi \Vert^{2}_{\mathcal{M}_{Yj}}$. For each $\xi \in \mathbb{R}^{2nq}$, there exist an $(x^{*},z^{*},v^{*},w^{*},y^{*}) \in \mathcal{M}_{Y}$ such that $\nabla V(\xi) = [(\varphi_{1}\bar{x}^{*} - \varphi_{2} \bar{z}^{j*})^{T}, (\varphi_{3} \bar{z}^{j*})^{T}]^{T}$.
\end{itemize}
\end{lem}

\begin{pf}
Obviously $(i)$ is true. It follows from $(i)$ that there exists a unique $(x^{*},z^{*},v^{*},w^{*},y^{*})$ such that $P_{\mathcal{M}_{Yj}}(\xi)=[(\varphi_{1}x^{*} - \varphi_{2}z^{j*})^{T}, (\varphi_{3}z^{j*})^{T}]^{T}$. By definition of $P_{\mathcal{M}_{Yj}}(\xi)$, $[(\varphi_{1}x^{*} - \varphi_{2}z^{j*})^{T}, (\varphi_{3}z^{j*})^{T}]^{T} \in \mathcal{M}_{Yj}$, which means that $(x^{*},z^{*},v^{*},w^{*},y^{*}) \in \mathcal{M}_{Y}$. Thus $(ii)$ is proved.

Note that $\Vert \xi \Vert^{2}_{\mathcal{M}_{Yj}} = \Vert \xi - P_{\mathcal{M}_{Yj}}(\xi) \Vert^{2}$. Then according to $(i)$ and $(ii)$, it shows that $\Vert \xi \Vert^{2}_{\mathcal{M}_{Yj}} = \Vert [((\varphi_{1}x - \varphi_{2}z^{j})-(\varphi_{1}x^{*} - \varphi_{2}z^{j*}))^{T}, (\varphi_{3}z^{j} - \varphi_{3}z^{j*})^{T}]^{T} \Vert^{2}=\Vert (\varphi_{1}x - \varphi_{1}x^{*})-(\varphi_{2}z^{j} - \varphi_{2}z^{j*}) \Vert^{2} + \Vert \varphi_{3}z^{j} - \varphi_{3}z^{j*} \Vert^{2}$. Hence $(iii)$ is proved.

Similarly to the analysis of $(iii)$, there holds that $\nabla V(\xi) = \frac{1}{2} \nabla \Vert \xi - P_{\mathcal{M}_{Yj}}(\xi) \Vert^{2} = \xi - P_{\mathcal{M}_{Yj}}(\xi)$. Then according to $(ii)$, it is shown that $\xi - P_{\mathcal{M}_{Yj}}(\xi) = [(\varphi_{1} \bar{x}^{*} - \varphi_{2} \bar{z}^{j*})^{T}, (\varphi_{3} \bar{z}^{j*})^{T}]^{T}$. This completes the proof of $(iv)$. $\hfill$ $\blacksquare$
\end{pf}

Then, the main theorem of this subsection is given below.

\begin{thm}
Consider algorithm \eqref{Algorithm 2}. Suppose Assumptions $\ref{A1}$-$\ref{A4}$ hold. If inequalities \eqref{IE} hold, then the trajectory of $x(t)$ converges, and $\lim\limits_{t \to \infty} x(t)$ is the solution of problem \eqref{Problem 2}. \label{Thm2}
\end{thm}

\begin{pf}
Define $\phi = col(x,z,v,w)$. The first-order system controlled by (\ref{Algorithm 2}) can be considered as
\begin{gather}\label{G123}
\dot{\phi} = g_{1}(\phi) + g_{2}(\phi,y) + g_{3}(y),
\end{gather}
where $g_{1}(\phi)=col(\dot{x}$,$\dot{z}$,$G_{1}$,$\dot{w})$, $G_{1} = - H^{-1}_{nq}(x - d) - \alpha L_{nq} v - w$, $g_{2}(\phi,y) = col(\textbf{0}_{nq},\textbf{0}_{(m-1)nq},G_{2},\textbf{0}_{nq})$, $G_{2} = (H^{-1}_{nq} \!\!-\!\! Y^{-1})\bar{x}^{*}$, $g_{3}(y) = col(\textbf{0}_{nq},\textbf{0}_{(m-1)nq},u,\textbf{0}_{nq})$, and $u(t) = (H^{-1}_{nq} \!\!-\!\! Y^{-1}(t))(x^{*} \!\!-\!\! d)$.

\textbf{i)} Firstly, with only the first part in \eqref{G123}, we consider the system
\begin{equation}
\dot{\phi} = g_{1}(\phi). \label{G1}
\end{equation}
From Theorem \ref{Thm1}, it is clear that under system \eqref{G1}, $(x(t), z(t), v(t), w(t))$ converges to the largest invariant set $\mathcal{M}$ in $ E = \lbrace (x,z,v,w) \vert x=x^{*}, v=v^{*},w=w^{*}, -\gamma z^{j} \in \partial F^{j}(x^{*})$ for $j \in \lbrace 1, \cdots, m-1\rbrace \rbrace$.

\textbf{ii)} Consider the system
\begin{equation}
\dot{\phi} = g_{1}(\phi) + g_{2}(\phi, y),
\label{G12}
\end{equation}
where $[(\phi^{*})^{T},(y^{*})^{T}]^{T}$ is an equilibrium of algorithm \eqref{Algorithm 2}, and $g_{2}(\phi, y)$ satisfies that $g_{2}(\phi^{*},y^{*}) = \textbf{0}$. 

From \eqref{G12} and the Lyapunov candidate $V_{Y}(x,z,v,w,y)=V(x,z,v,w)+V_{4}(y)$, where $V_{4}(y) = \frac{1}{2} \Vert \bar{y}^{*} \Vert^{2}$ and $\bar{y}^{*} \triangleq y - y^{*}$, it yields that
\begin{gather}\label{DVG12}
\begin{split}
& \dot{V}_{Y}(x,z,v,w,y)\\
\leq & \! - \epsilon_{1} \Vert \dot{x} \Vert^{2} \!\!-\! \epsilon_{2} \sum_{j=1}^{m-1} \Vert \dot{z}^{j} \Vert^{2} \!\! - \! \epsilon_{3}\Vert \bar{x}^{*} \Vert^{2} \!\! - \! \epsilon_{4} \Vert \bar{v}^{*} \! \Vert^{2}\\
& - \epsilon_{5} (\bar{w}^{*})^{T} \! H_{nq}\bar{w}^{*} - \frac{1}{2}(\bar{y}^{*})^{T}(L_{nn}+L_{nn}^{T})\bar{y}^{*} \\
& + DV_{Y},\\
\end{split}
\end{gather}
where 
\begin{gather}\label{DVG2}
\begin{split}
& DV_{Y} =  \frac{\partial V_{3}(v,w)}{\partial v} G_{2} \\
= & (\eta+1) (\bar{v}^{*})^{T}Qx^{-*} + (\bar{w}^{*})^{T}Qx^{-*} \\
\leq &  \zeta_{1} (\bar{v}^{*})^{T}Qv^{-*} + \zeta_{2} (\bar{x}^{*})^{T}Qx^{-*} + \zeta_{3} (\bar{w}^{*})^{T}Qw^{-*}\\
\leq & \rho(t) \left[ \zeta_{1} \Vert \bar{v}^{*} \Vert^{2} + \zeta_{2} \Vert \bar{x}^{*} \Vert^{2} + \zeta_{3} \Vert \bar{w}^{*} \Vert^{2} \right]\\
\end{split}
\end{gather}
and $Q = I_{nq} - H_{nq}Y^{-1}$, $\zeta_{1} = \frac{\eta+1}{2}$, $\zeta_{2} = \frac{\eta}{2}+1$, $\zeta_{3} = \frac{1}{2}$, $\rho(t) = \max_{i\in I}\vert 1 - h_{i}(y_{i}^{i}(t))^{-1} \vert$.

Since $y(t) = e^{-L_{nn}t}y(0)$ and $y(0)=[I_{n}^{1},\cdots,I_{n}^{n}]^{T}$ from \eqref{Algorithm 2}, it is shown that $\lim\limits_{t \to \infty} y(t) = \textbf{1}(h^{T} \otimes I_{q}) y(0) = \textbf{1}_{n} \otimes h$. Therefore, $y^{*} = \textbf{1}_{n} \otimes h$. Then according to Lemma 2.6 in \cite{RWB}, $y(t)$ is exponentially convergent to $\textbf{1}_{n}\otimes h$, and $y_{i}^{i}(t) > 0$ for all $i \in \lbrace 1,\cdots, n\rbrace$ and $t > 0$. As the result, $\rho(t)$ and $u(t)$ are both exponentially convergent to zero.

With \eqref{DVG12} and \eqref{DVG2}, it is followed that
\begin{gather}
\begin{split}
& \dot{V}_{Y}(x,z,v,w,y)\\
\leq & - \epsilon_{1} \Vert \dot{x} \Vert^{2} - \epsilon_{2} \sum_{j=1}^{m-1} \Vert \dot{z}^{j} \Vert^{2} - l_{1}\Vert \bar{x}^{*} \Vert^{2}\\
& - l_{2} \Vert \bar{v}^{*} \Vert^{2} - l_{3} (\bar{w}^{*})^{T}H_{nq}\bar{w}^{*},\\
\end{split}
\end{gather}
where $l_{1}=\epsilon_{3}-\rho(t)\zeta_{2}$, $l_{2}=\epsilon_{4}-\rho(t)\zeta_{1}$, $l_{3}=\epsilon_{5}-\rho(t)\frac{\zeta_{3}}{h^{*}}$.

Since $\rho(t) \to 0$ when $t \to \infty$, there exists $T_{0} > 0$ that when $t > T_{0}$, $ l_{1} \geq \frac{1}{2}\epsilon_{3}, l_{2} \geq \frac{1}{2}\epsilon_{4} , l_{3} \geq \frac{1}{2}\epsilon_{5}$. Therefore $\dot{V}_{Y}(x,z,v,w,y) \leq 0$ when $t > T_{0}$.

When $t \leq T_{0}$, since $0 < \rho(t) < 1$, 
\begin{gather}
\begin{split}
& \dot{V}_{Y}(x,z,v,w,y) \\
\leq &  \zeta_{2}\Vert \bar{x}^{*} \Vert^{2} + \zeta_{1} \Vert \bar{v}^{*} \Vert^{2} + \zeta_{3} \Vert \bar{w}^{*}\Vert^{2}\\
\leq &  \zeta_{2}(\Vert \bar{v}^{*} \Vert^{2}+\Vert \bar{w}^{*}\Vert^{2}) + \zeta_{2} \Vert \bar{x}^{*} \Vert^{2} \\
\leq & \zeta_{2}(\Vert \bar{v}^{*} \Vert^{2}+\Vert \bar{w}^{*}\Vert^{2}) + \iota_{1} V_{1}(x,z), \label{DVT}
\end{split}
\end{gather}
where $\iota_{1} = \frac{\eta+2}{(\eta+1)[1-(m-1)\gamma]}$.

Note that
\begin{gather}\label{V3T}
\begin{split}
& V_{3}(v,w) \\
\geq & \frac{\eta+1}{2h_{max}}\Vert \bar{v}^{*} \Vert^{2} + \frac{1}{2h_{max}} \Vert \bar{w}^{*} \Vert^{2} + \frac{1}{h_{max}}(\bar{v}^{*})^{T}\bar{w}^{*} \\
\geq &  \iota_{2}\Vert \bar{v}^{*} \Vert^{2} + \iota_{3} \Vert \bar{w}^{*} \Vert^{2} \\
\geq &  \min \lbrace \iota_{2}, \iota_{3} \rbrace (\Vert \bar{v}^{*} \Vert^{2} + \Vert \bar{w}^{*} \Vert^{2}), 
\end{split}
\end{gather}
where $h_{max} = max_{i \in I} \lbrace h_{1}, \cdots, h_{n} \rbrace$, $1< r < \eta + 1$, $\iota_{2} = \frac{\eta+1}{2h_{max}} - \frac{1}{2h_{max} r}$, $\iota_{3} = \frac{1}{2h_{max}} - \frac{r}{2}$.

From \eqref{DVT} and \eqref{V3T}, it is shown that
\begin{gather}\label{DVTK}
\begin{split}
& \dot{V}_{Y}(x,z,v,w,y) \\
\leq &  \frac{\zeta_{2}}{\min \lbrace \iota_{2}, \iota_{3} \rbrace} V_{3}(x,z,v,w) + \iota_{1} V_{1}(x,z) \\
\leq & \kappa_{1} V_{Y}(x,z,v,w,y), 
\end{split}
\end{gather}
where $\kappa_{1} = \max \lbrace \frac{\zeta_{2}}{\min \lbrace \iota_{2}, \iota_{3} \rbrace}, \iota_{1} \rbrace$. 

According to \eqref{DVTK}, when $t=T_{0}$,
\begin{equation}
V_{Y}(T_{0}) \leq e^{\kappa T_{0}} V_{Y}(0).
\label{VG12}
\end{equation}

To sum up, $V_{Y}(x,z,v,w,y)$ is proper and $\dot{V}_{Y}(x,z,v,w,y)$ $\leq 0$ when $t > T_{0}$. Then according to the LaSalle invariant principle and Lemma \ref{Aequilibrium}, system \eqref{G12} converges to the largest invariant set $\mathcal{M}_{Y}$ in $ E = \lbrace (x,z,v,w,y) \vert x=x^{*}, v=v^{*},w=w^{*}, y=y^{*}, -\gamma z^{j} \in \partial F^{j}(x^{*})$ for $j \in \lbrace 1, \cdots, m-1\rbrace \rbrace$, which also means that system \eqref{G12} is GAS to $\mathcal{M}_{Y}$.

\textbf{iii)} Now consider the complete system \eqref{G123}. Clearly $\mathcal{M}_{Y}$ is a closed, 0-invariant set for system \eqref{G123}. Similar to \eqref{VG12}, there exist $\kappa_{2} > 0$ and $\nu_{u} >0$ such that
\begin{gather}\label{VG123}
\begin{split}
& \forall t \leq T_{0}:\quad \dot{V}_{Y}(t) \leq \kappa_{2} V_{Y}(t) + \nu_{u} \Vert u(t) \Vert^{2},\\
& \forall t > T_{0}:\quad \dot{V}_{Y}(t) \leq  \nu_{u} \Vert u(t) \Vert^{2}.
\end{split}
\end{gather}
As a result, $V_{Y}(t)$ is bounded for all $t < +\infty$. Since $V_{Y}(t)$ is proper, system \eqref{G123} is forward complete. Moreover, note that $T_{0} \to 0$ and $\Vert u(t) \Vert \to 0$ when $y \to y^{*}$. Therefore, according to \eqref{VG123}, each $(x^{*},z^{*},v^{*},w^{*},y^{*}) \in \mathcal{M}_{Y}$ is Lyapunov stable.

Then, we define an iISS-Lyapunov candidate $V_{\mathcal{M}_{Y}}(x,z,$ $v,w,y) = V_{1\mathcal{M}_{Y}}(x,z)+V_{2}(x)+V_{3}(v,w)+V_{4}(y)$ with respect to $\mathcal{M}_{Y}$, where
\begin{gather}
V_{\!1\!\mathcal{M}_{\!Y}}\!(\!x,\!z\!) \!\!=\!\! \frac{\eta \!+\!  1}{2} \!\!\sum_{j=1}^{m-1}\!\!  \Vert [ (\varphi_{1} x \!-\! \varphi_{2} z^{j})^{\!T}\!\!, (\varphi_{3} z^{j})^{\!T}]^{\!T}\!\Vert^{2}_{\!\mathcal{M}_{\!Yj}},
\end{gather}
and $\varphi_{1}=(\frac{1}{m-1})^{\frac{1}{2}}$, $\varphi_{2}=\gamma (m-1)^{\frac{1}{2}}$, $\varphi_{3}=[\gamma (1 - \gamma(m-1))]^{\frac{1}{2}}$, $\mathcal{M}_{Yj} \!\! \triangleq \!\! \lbrace  [ (\varphi_{1} x^{*} - \varphi_{2} z^{j*})^{T} , (\varphi_{3} z^{j*})^{T} ]^{T}, (x^{*},z^{j*},v^{*}$, $w^{*},y^{*}) \in \mathcal{M}_{Y} \rbrace$ for $j \in \lbrace 1, \cdots, m-1 \rbrace$.

According to Lemma \ref{x-z} and proof of Theorem \ref{Thm1}, when $t > T_{0}$, it follows that
\begin{gather}
\begin{split}
& \dot{V}_{\mathcal{M}_{Y}}(x,z,v,w,y) \\
\leq & - l_{2} \Vert \bar{v}^{*} \Vert^{2} - l_{3} (\bar{w}^{*})^{T}H_{nq}\bar{w}^{*} \\
&  + \left[ (\bar{v}^{*} + \bar{w}^{*})^{T}H_{nq} + \eta(\bar{v}^{*})^{T}H_{nq} \right] u(t) \\
\leq & - \iota_{4} \Vert \bar{v}^{*} \Vert^{2} - \iota_{5}(\bar{w}^{*})^{T}H_{nq}\bar{w}^{*} + \iota_{6} \Vert u(t) \Vert^{2} \\
\leq & \iota_{6} \Vert u(t) \Vert^{2}, \label{0P}
\end{split}
\end{gather}
where $\iota_{4} = \frac{1}{2}\epsilon_{3} - \frac{(\eta+1)\tau_{1}}{2h^{*}}$, $\iota_{5} = \frac{1}{2}\epsilon_{4} - \frac{\tau_{2}}{2}$ and $\iota_{6} = \frac{\eta+1}{2\tau_{1}} + \frac{1}{2\tau_{2}}$. Note that $\iota_{4} > 0$ and $\iota_{5} > 0$ always hold, since $\tau_{1}$ and $\tau_{2}$ can be chosen arbitrarily small.

From \eqref{0P} and Definition $\ref{D2}$, it is clear that system \eqref{G123} is ZOD with respect to $\mathcal{M}_{Y}$ when $t > T_{0}$. Remind that $u(t)$ is exponentially convergent to zero. Since system \eqref{G123} is 0-GAS with respect to $\mathcal{M}_{Y}$, then according to Lemma \ref{Diiss}, $\phi_{Y}(t,u(t))$ converges to $\mathcal{M}_{Y}$. Note that each $(x^{*},z^{*},v^{*},w^{*},y^{*}) \in \mathcal{M}_{Y}$ is Lyapunov stable, then according to Lemma \ref{Aequilibrium}, system \eqref{G123} converges to one of its equilibria in $\mathcal{M}_{Y}$. Similar to the analysis in proof of Theorem $\ref{Thm1}$, it is clear that $(\textbf{1}_{n}\otimes I_{q})^{T}H_{nq}w^{*} = \textbf{0}_{q}$. Then according to Lemma $\ref{KA2}$, $x^{*}$ is the solution of problem (\ref{Problem 2}). This completes the proof. $\hfill$ $\blacksquare$
\end{pf}

\begin{rem}
For $\xi = [(\varphi_{1} x - \varphi_{2} z^{j})^{T} , (\varphi_{3} z^{j})^{T} ]^{T}$ with $j \in \lbrace 1, \cdots, m-1\rbrace$, let $P_{\mathcal{M}_{Yj}}(\xi)=[(\varphi_{1} \hat{x}^{*} - \varphi_{2} \hat{z}^{j*})^{T}, (\varphi_{3} \hat{z}^{j*})^{T} ]^{T}$, $P_{\mathcal{M}_{Yx}}(x) = \tilde{x}^{*}$, and $P_{\mathcal{M}_{Yz^{j}}}(z^{j}) = \tilde{z}^{j*}$, where $\mathcal{M}_{Yx} \triangleq \lbrace x^{*} \vert (x^{*},z^{*},v^{*},w^{*},y^{*}) \in \mathcal{M}_{Y} \rbrace$ and $\mathcal{M}_{Yz^{j}} \triangleq \lbrace z^{j*} \vert (x^{*},z^{*},v^{*},w^{*},y^{*}) \in \mathcal{M}_{Y} \rbrace$. Hence $\hat{x}^{*} = \tilde{x}^{*} = x^{*}$ and usually $\hat{z}^{j*} \neq \tilde{z}^{j*}$. While it is true that $\hat{z}^{j*} \in \mathcal{M}_{Yz^{j}}$, hence $\dot{V}_{\mathcal{M}_{Y}}(x,z,v,w,y)$ can be deduced based on the analysis of $\dot{V}_{Y}(x,z,v,w,y)$ in proof of Theorem \ref{Thm1}. $\hfill$ $\blacklozenge$
\end{rem}

\begin{rem}
In proof of Theorem $\ref{Thm2}$, the first-order system controlled by \eqref{Algorithm 2} had been separated to three parts. Since the existence of estimation error between $y$ and $h$, the Lyapunov function $V_{Y}(x,z,v,w,y)$ of system \eqref{G12} may increase before $T_{0}$. Then we proved that the Lyapunov function $V_{Y}(x,z,v,w,y)$ of system \eqref{G12} is bounded when $t \leq T_{0}$ and $\dot{V}_{Y}(x,z,v,w,y) \leq 0$ when $t > T_{0}$. Finally, with the help of iISS theory with respect to set, it is proved that system \eqref{G123} is asymptotically convergent to its equilibria in $\mathcal{M}_{Y}$, which provides new ideas about stability analysis of asymptotically convergent system with exponentially convergent inputs. $\hfill$ $\blacklozenge$
\end{rem}

\section{Simulations}
In this section, simulations are performed to validate the proposed algorithm \eqref{Algorithm 2}. Consider the fused LASSO problem with four agents moving in a 2-D space with first-order dynamics \eqref{system} as
\begin{gather}
\min_{x \in \mathbb{R}^{8}} F(x), \ \ s.t. \sum_{i=1}^{4} x_{i} = \sum_{i=1}^{4} d_{i},
\end{gather}
where $x_{i} = [ x^{1}_{i}, x^{2}_{i} ]^{T} \in \mathbb{R}^{2}, i \in \lbrace 1, 2, 3, 4 \rbrace$, $F(x) = \sum_{j=0}^{3} f^{j}(x) = 2\Vert x - s \Vert^{2} + \iota(x) + \Vert x - p \Vert_{1} + \Vert Dx \Vert_{1}$, $\iota(x) = \begin{cases} 
0,  & \mbox{if } x \in \Omega\\
\infty, & \mbox{if } x \notin \Omega
\end{cases}$, and 
\begin{equation}
D = 
\begin{bmatrix}
1  & -1  &      &      &        &     &    &\\
   &     &  1   & -1   &        &     &    &\\
   &     &      &      & \cdots & \cdots &  &  \\
   &     &      &      &        &     & 1 & -1\\
\end{bmatrix} \in \mathbb{R}^{8 \times 8}.
\end{equation} 
The local cost function $f_{i}(x_{i})$ for agent $i$ is consisted by
\begin{equation}
\begin{array}{l}
f^{0}_{i}(x_{i}) = 2\Vert x_{i} - s_{i}\Vert^{2},\\
f^{1}_{i}(x_{i}) = \Vert x_{i} - p_{i} \Vert_{1},\\
f^{2}_{i}(x_{i}) = \Vert x_{i}^{1} - x_{i}^{2} \Vert_{1},\\
f^{3}_{i}(x_{i}) =
\begin{cases} 
0,  & \mbox{if } x_{i} \in \Omega_{i}\\
\infty, & \mbox{if } x_{i} \notin \Omega_{i}
\end{cases},\\
\end{array}
\label{Simulation_f}
\end{equation}
where $s_{i} = [s^{1}_{i}, s^{2}_{i}]^{T} = [i - 2.5, 0]^{T}$, $p_{i} = [p^{1}_{i}, p^{2}_{i}]^{T} = [0, i - 2.5]^{T}$ and $\Omega_{i} = \lbrace \delta \in \mathbb{R}^{2} \vert \Vert \delta - x_{i}(0) \Vert^{2} \leq 64\rbrace$. Then $f^{0}_{i}(x_{i})$, $f^{1}_{i}(x_{i})$, $f^{2}_{i}(x_{i})$ and $f^{3}_{i}(x_{i})$ represent respectively the quadratic objective, the $l_{1}$ penalty with an anchor $p_{i}$, another $l_{1}$ penalty associated with the matrix $D$, and the indicator function of the constraint set $x_{i} \in \Omega_{i}$ for each agent $i$. Resource allocation conditions are described as $d_{1}=[2,-1]^{T}$, $d_{2}=[-1,1]^{T}$, $d_{3}=[-1,-1]^{T}$ and $d_{4}=[2,2]^{T}$.

Based on ($\ref{Simulation_f}$), the gradient of $f^{0}_{i}$ and proximal operators of $f^{1}_{i}$, $f^{2}_{i}$ and $f^{3}_{i}$ for agent $i$ are shown as
\begin{gather}
\begin{split}
\nabla f^{0}_{i}(x_{i}) = & [ 4(x^{1}_{i} - s^{1}_{i}), 4(x^{2}_{i} - s^{2}_{i}) ]^{T}, \\
prox_{f^{1}_{i}}[\eta_{1}] = & [\phi(\eta_{1}^{1}, p_{i}^{1}), \phi(\eta_{1}^{2}, p_{i}^{2})]^{T}, \\
prox_{f^{2}_{i}}[\eta_{2}] = & [\phi(\eta_{2}^{1}, \eta_{2}^{2}), \phi(\eta_{2}^{2}, \eta_{2}^{1})]^{T}, \\
prox_{f^{3}_{i}}[\eta_{3}] = & \arg \min_{\delta \in \Omega_{i} } \Vert \delta - \eta_{3} \Vert^{2}, \\
\end{split}
\end{gather}
where $\eta_{j} \in \mathbb{R}^{2}$, $j \in \lbrace 1, 2 ,3 \rbrace$. For $\xi_{1} \in \mathbb{R}$ and $ \xi_{2} \in \mathbb{R}$, the function $\phi(\xi_{1},\xi_{2})$ is defined as follows
\begin{eqnarray}
\phi(\xi_{1},\xi_{2}) & = &
\begin{cases} 
\xi_{1}-1,  & \mbox{if } \xi_{1} > \xi_{2}+1 \\
\xi_{2},  & \mbox{if } \xi_{2}-1 \leq \xi_{1} \leq \xi_{2}+1 \\
\xi_{1}+1,  & \mbox{if } \xi_{1} <  \xi_{2}-1
\end{cases}.
\end{eqnarray}

Note that the proximal operator of $f^{1}_{i}(x_{i})+f^{2}_{i}(x_{i})+f^{3}_{i}(x_{i})$, e.i., $ prox_{(f^{1}_{i}+f^{2}_{i} + f^{3}_{i})}[\eta_{4}] = \arg \min_{\delta \in \Omega_{i}} \lbrace \Vert \delta \! - \! p_{i} \Vert_{1} + \Vert \delta^{1} \! - \! \delta^{2} \Vert_{1} + \frac{1}{2} \Vert \delta - \eta_{4} \Vert^{2} \rbrace$ is not proximable, where $\eta_{4} \in \mathbb{R}^{2}$. Hence proximal algorithms \cite{DCFO1}-\cite{OFW} may not fit for this problem.

The Laplacian matrix of weight-unbalanced directed graph $\mathcal{G}$ is given as 
\begin{equation}
L_{4} = 
\begin{bmatrix}
 1   &  0   &  0  & -1  \\
-1   &  2   & -1  &  0  \\
 0   & -1   &  1  &  0  \\
 0   &  0   & -1  &  1   \\
\end{bmatrix}
\end{equation}

We set $\alpha = 5$ and $\gamma = 0.2$ as coefficients in algorithm \eqref{Algorithm 2}. Initial positions of agents 1, 2, 3, and 4 are set as $x_{1}(0) = [-4, 5.5]^{T}$, $x_{2}(0) = [6, 5]^{T}$, $x_{3}(0) = [5, -3.5]^{T}$, and $x_{4}(0) = [-5, -5]^{T}$. We set initial values for Lagrangian multipliers $v_{i}$ and auxiliary variables $z_{i}^{1}, z_{i}^{2}, w_{i}$ for $i \in \lbrace 1, 2, 3, 4 \rbrace$ as zeros.

\begin{figure}
\centering
\subfigure{
\includegraphics[width=0.20\textwidth]{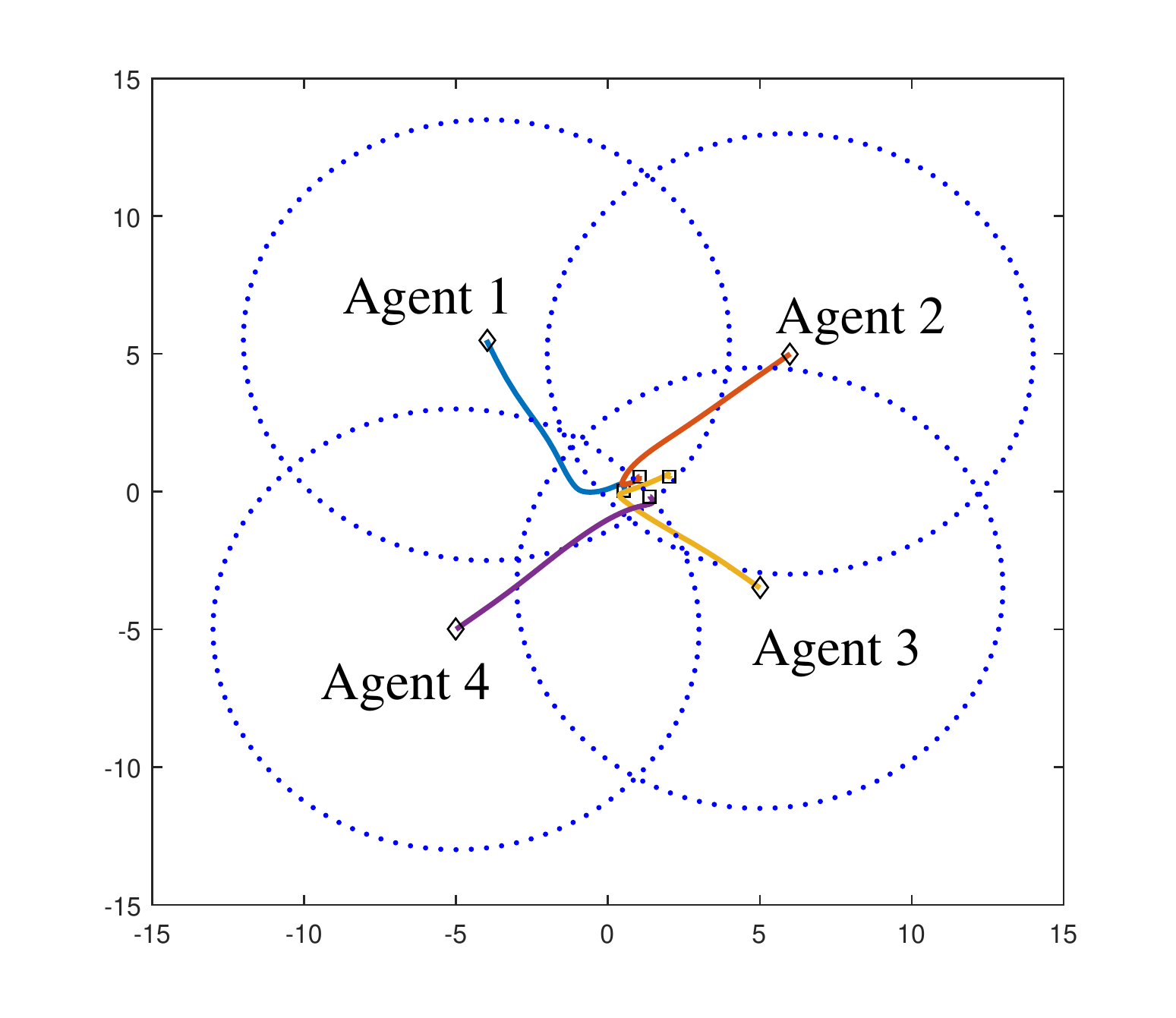}}
\subfigure{
\includegraphics[width=0.20\textwidth]{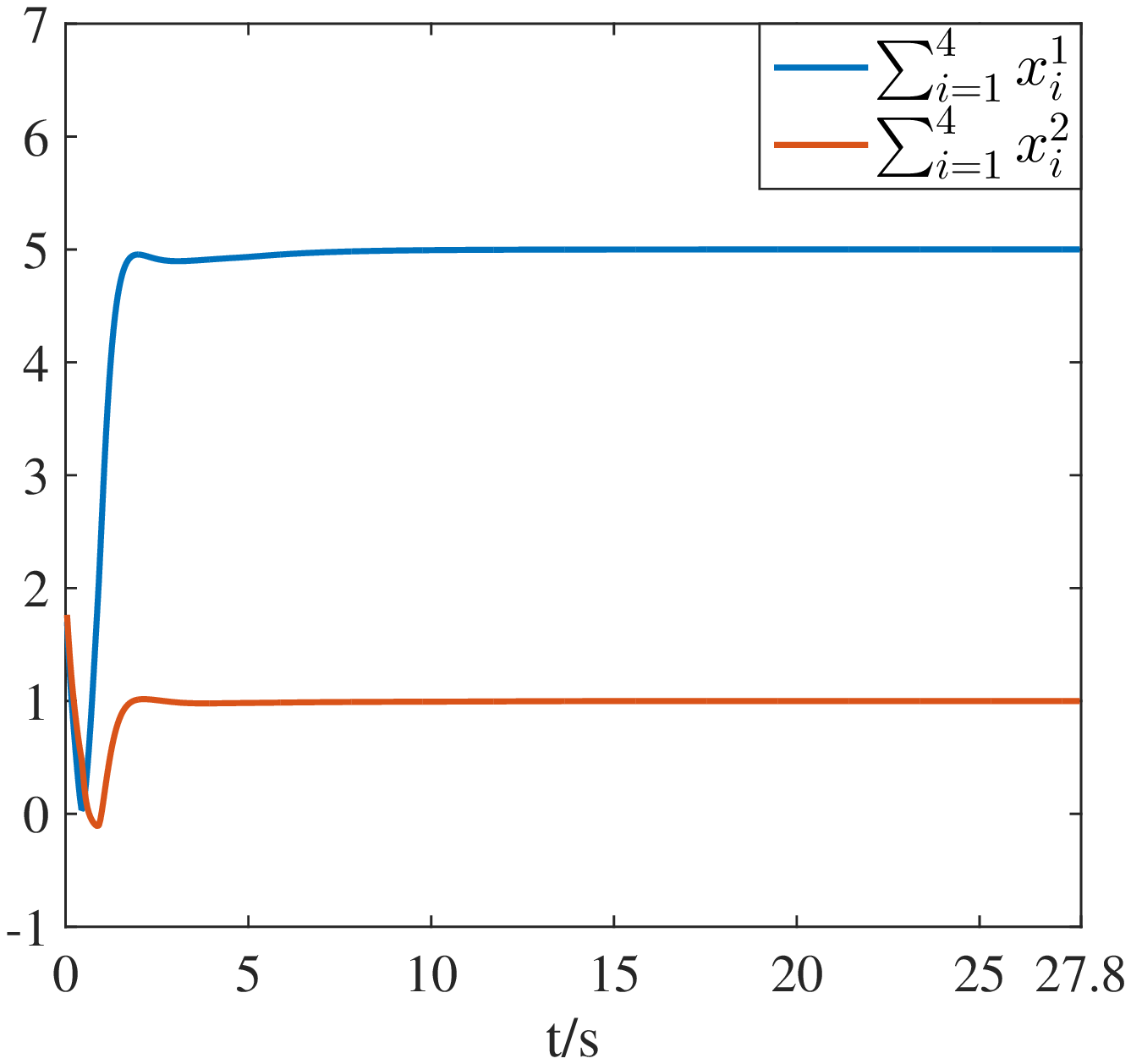}}
\caption{Motions of system \eqref{system} in a 2-D space and trajectories of $\sum_{i=1}^{4} x^{j}_{i}$ for $j \in \lbrace 1, 2 \rbrace$ with algorithm \eqref{Algorithm 2}}
\label{Fig.1}
\end{figure}

\begin{figure}
\centering
\subfigure{
\includegraphics[width=0.2\textwidth]{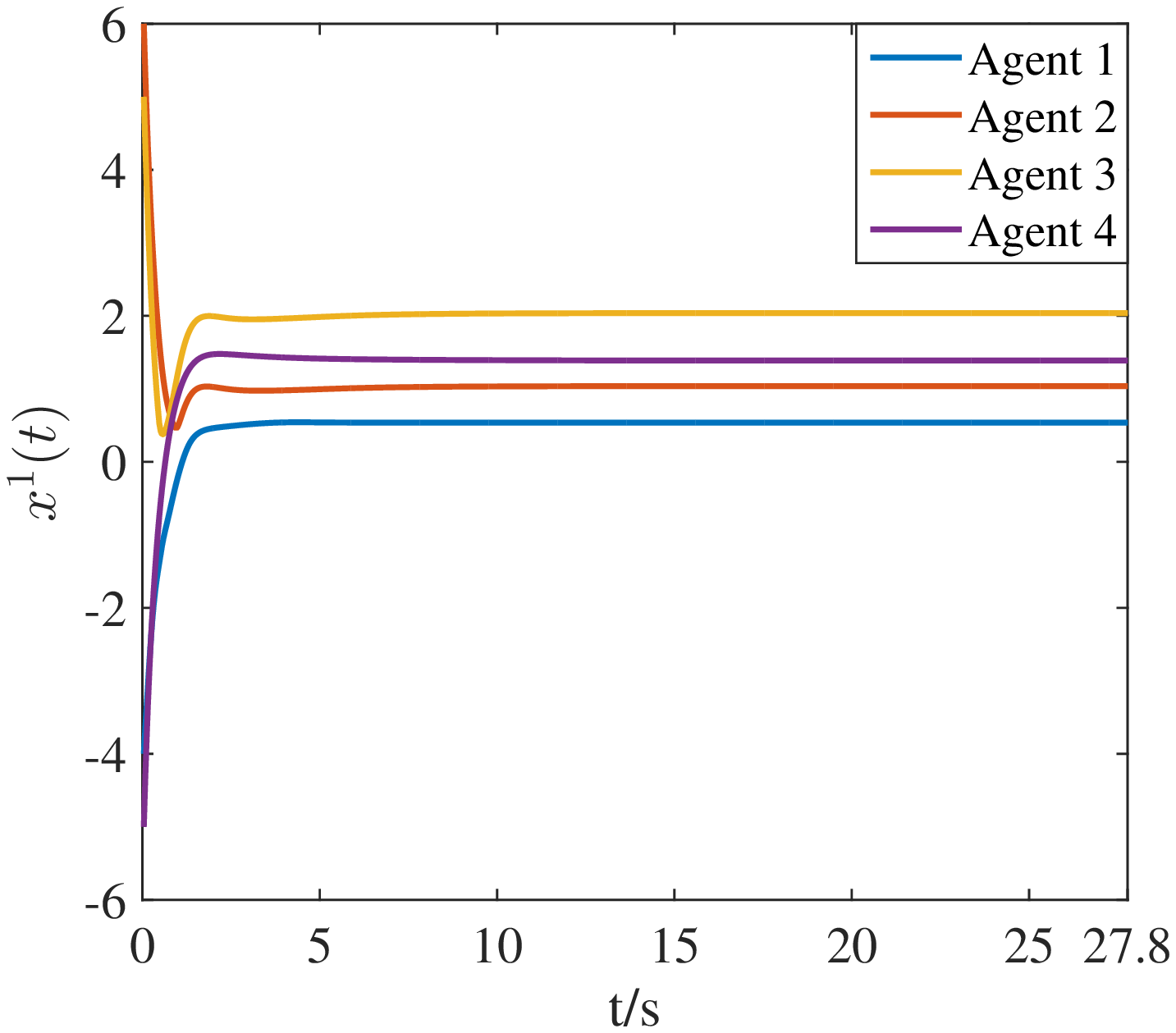}}
\subfigure{
\includegraphics[width=0.2\textwidth]{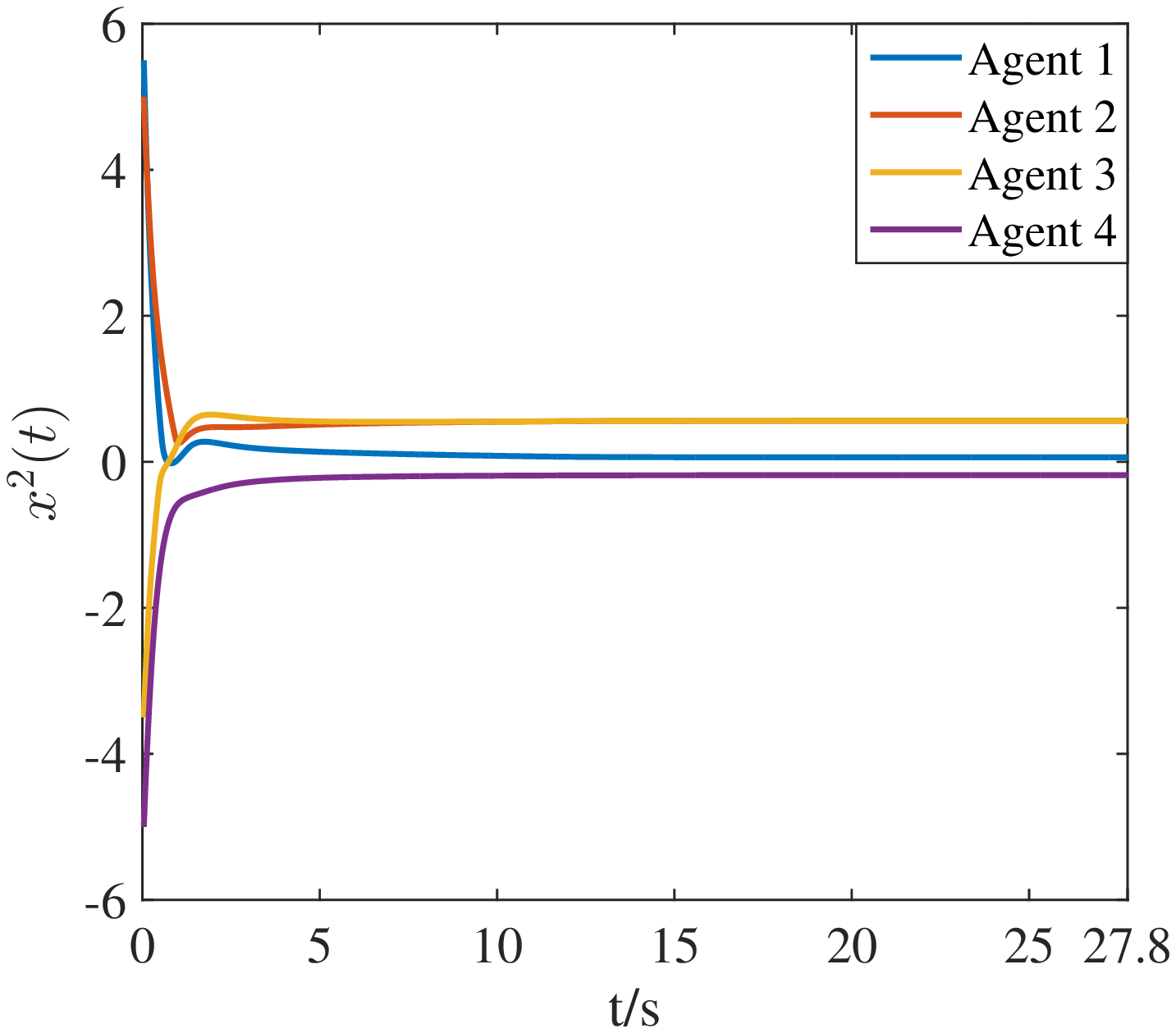}}
\caption{Trajectories of $x_{i}(t)$ for $i \in \lbrace 1, 2, 3, 4 \rbrace$ with algorithm ($\ref{Algorithm 2}$)}
\label{Fig.2}
\end{figure}

\begin{figure}[h]
\centering
\includegraphics[width=0.2\textwidth]{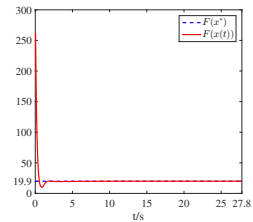}
\caption{The trajectory of $F(x)$ with algorithm ($\ref{Algorithm 2}$)}
\label{Fig.3}
\end{figure}

\begin{figure}
\centering
\subfigure{
\includegraphics[width=0.2\textwidth]{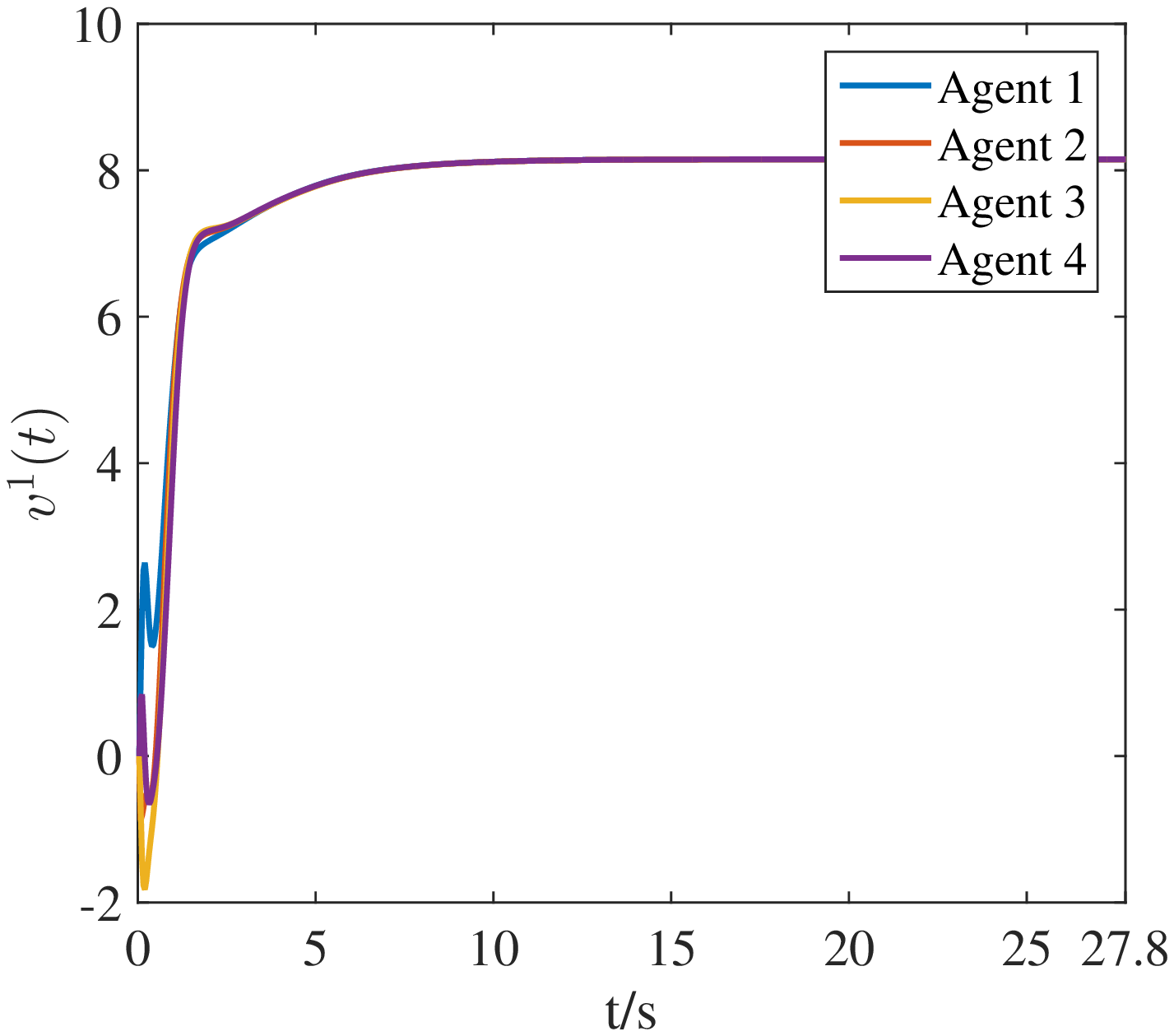}}
\subfigure{
\includegraphics[width=0.2\textwidth]{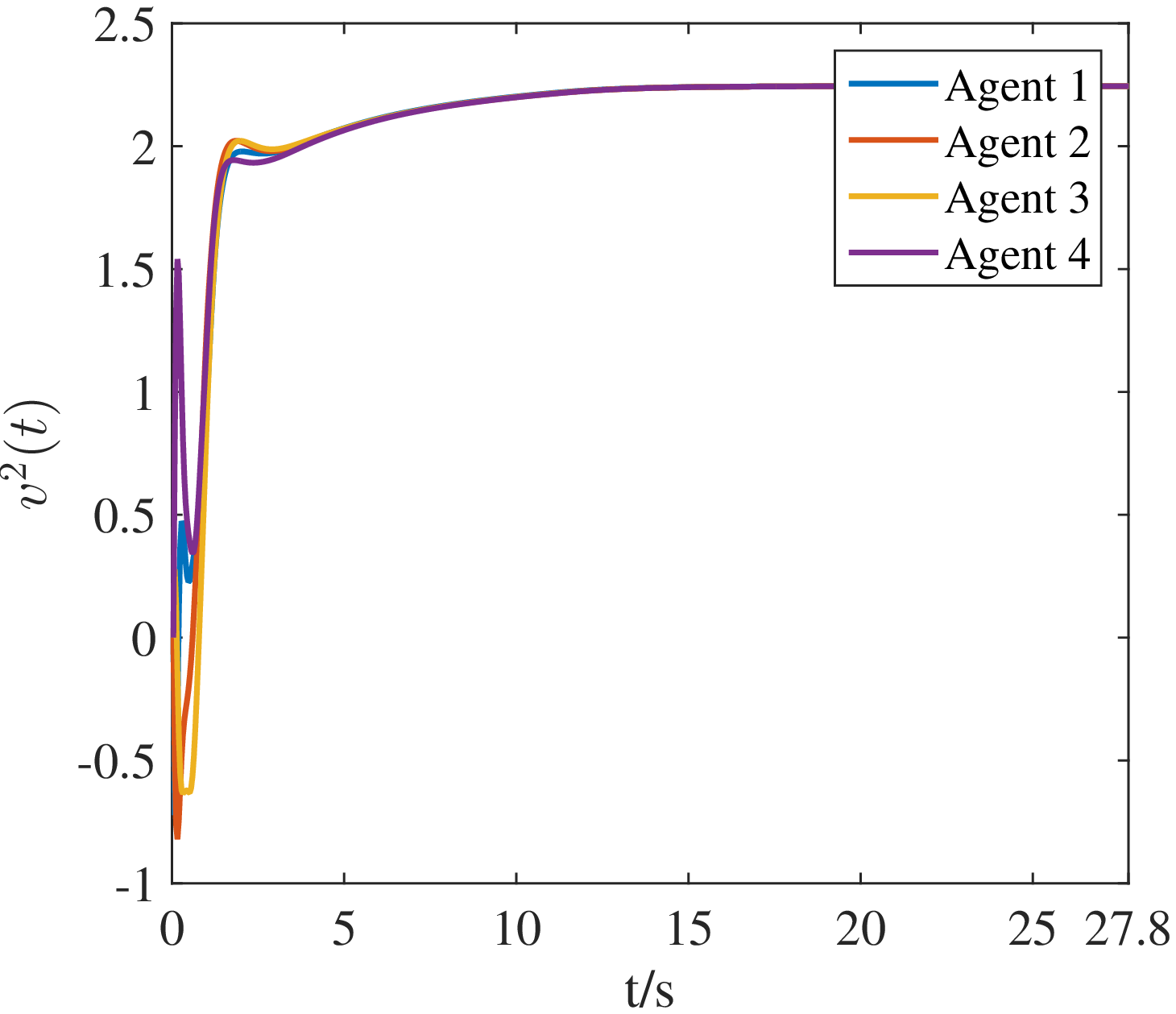}}
\caption{Trajectories of $v_{i}(t)$ for $i \in \lbrace 1, 2, 3, 4 \rbrace$ with algorithm ($\ref{Algorithm 2}$)}
\label{Fig.4}
\end{figure}

Motions of system \eqref{system} versus time and trajectories of $\sum_{i=1}^{4} x^{1}_{i}$ and $\sum_{i=1}^{4} x^{2}_{i}$ with algorithm ($\ref{Algorithm 2}$) are shown in Fig.$\ref{Fig.1}$, which show that resource allocation conditions $\sum_{i=1}^{4} x^{1*}_{i} = \sum_{i=1}^{4} d_{i}^{1} = 5$ and $\sum_{i=1}^{4} x^{2*}_{i} = \sum_{i=1}^{4} d_{i}^{2} = 1$ are satisfied. Fig.$\ref{Fig.2}$ gives trajectories of $x_{i}(t)$ for $i \in \lbrace 1, 2, 3, 4 \rbrace$. Fig.$\ref{Fig.3}$ shows the trajectory of $F(x)$, which proves that the global cost function is minimized. It can be seen from Fig.$\ref{Fig.1}$-Fig.$\ref{Fig.3}$ that all agents converge to the optimal solution which minimizes the global cost function and satisfies resource allocation conditions. Fig.$\ref{Fig.4}$ - Fig.$\ref{Fig.7}$ show trajectories of Lagrange multipliers $v_{i}(t)$ and auxiliary variables $z_{i}^{1}, z_{i}^{2}, w_{i}$ for $i \in \lbrace 1, 2, 3, 4 \rbrace$ respectively, which also verify the boundedness of system \eqref{system} steered by algorithm \eqref{Algorithm 2}. 


\begin{figure}
\centering
\subfigure{
\includegraphics[width=0.20\textwidth]{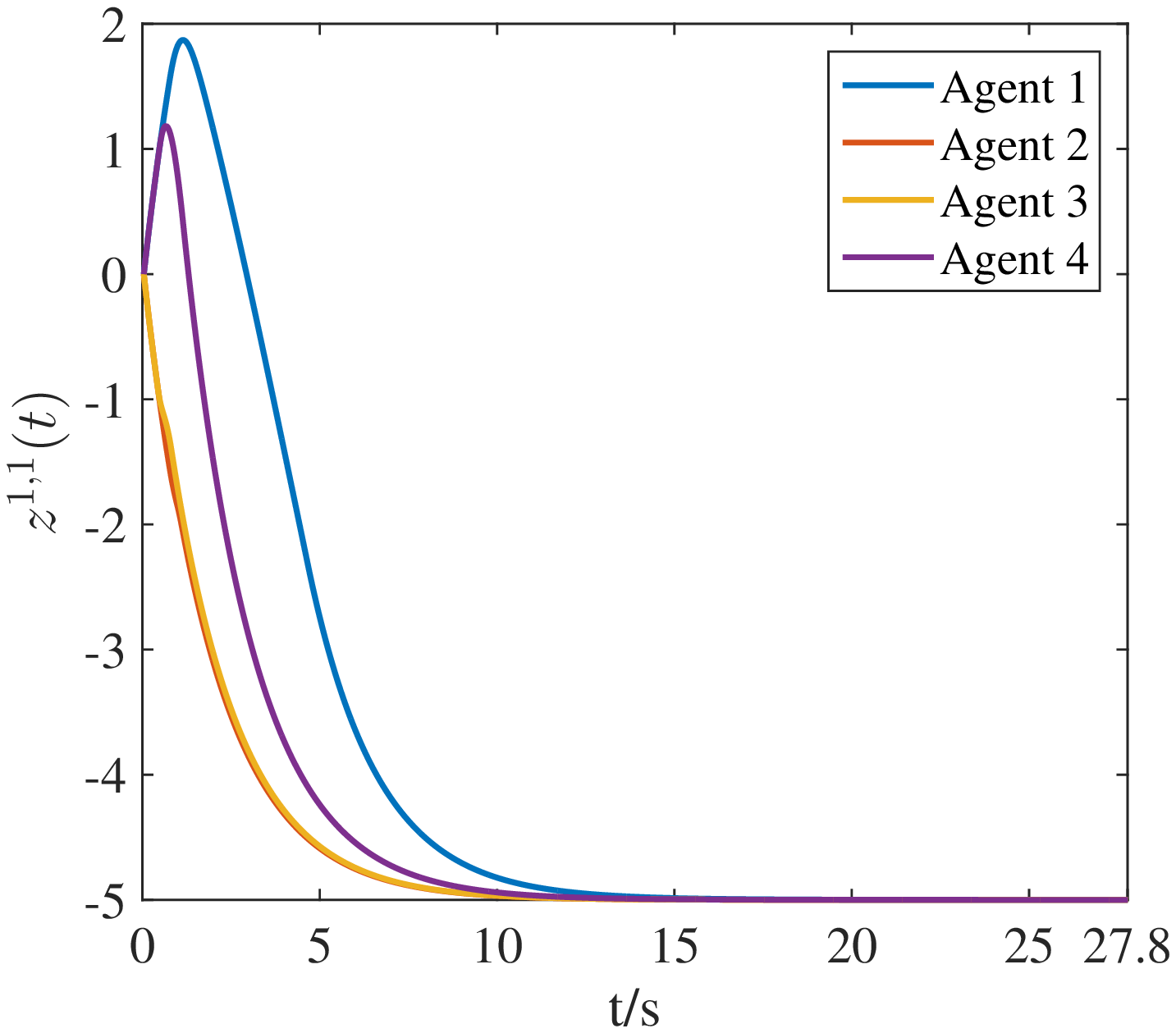}}
\subfigure{
\includegraphics[width=0.20\textwidth]{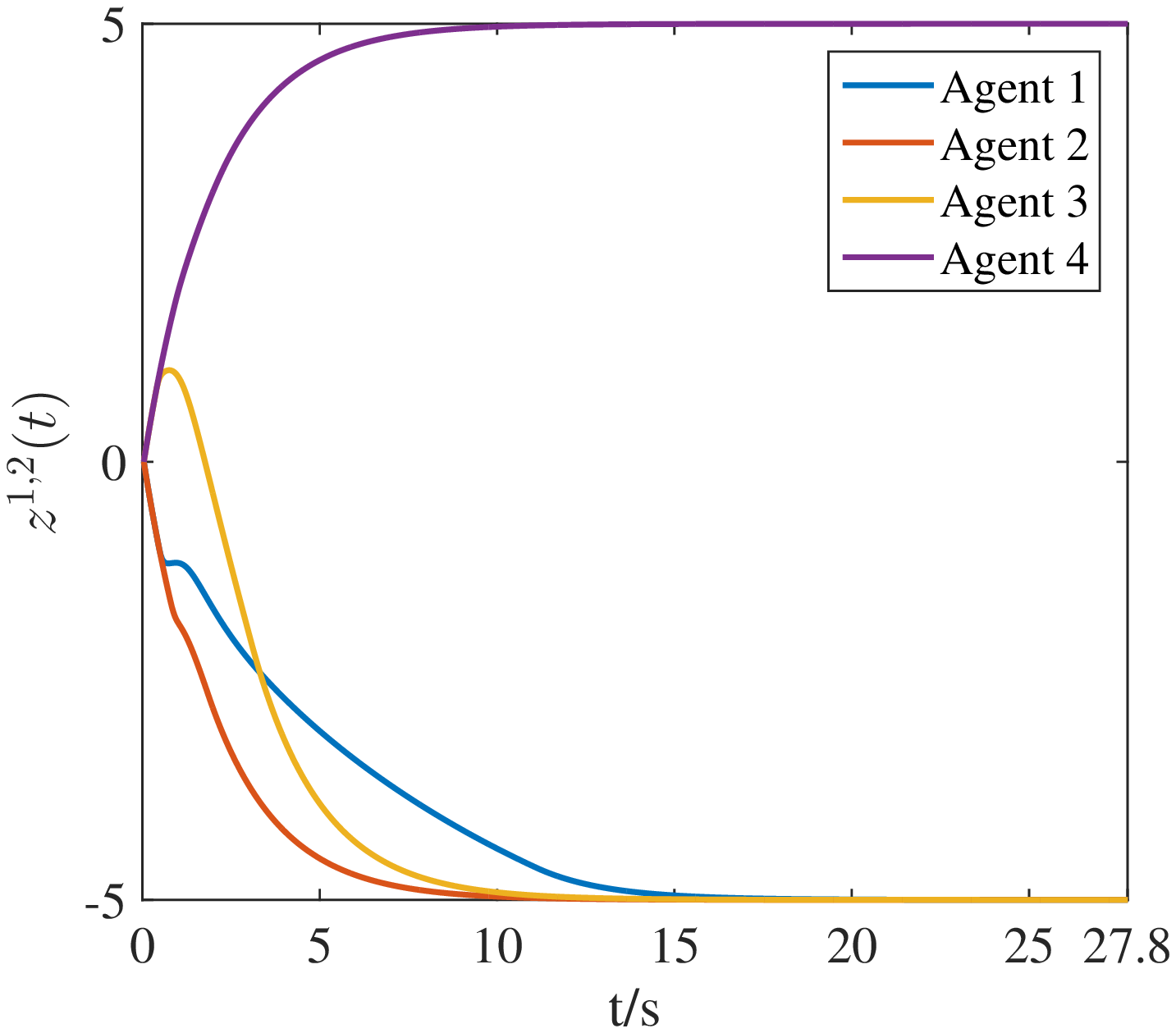}}
\caption{Trajectories of $z_{i}^{1}(t)$ for $i \in \lbrace 1, 2, 3, 4 \rbrace$ with algorithm ($\ref{Algorithm 2}$)}
\label{Fig.5}
\end{figure}

\begin{figure}
\centering
\subfigure{
\includegraphics[width=0.20\textwidth]{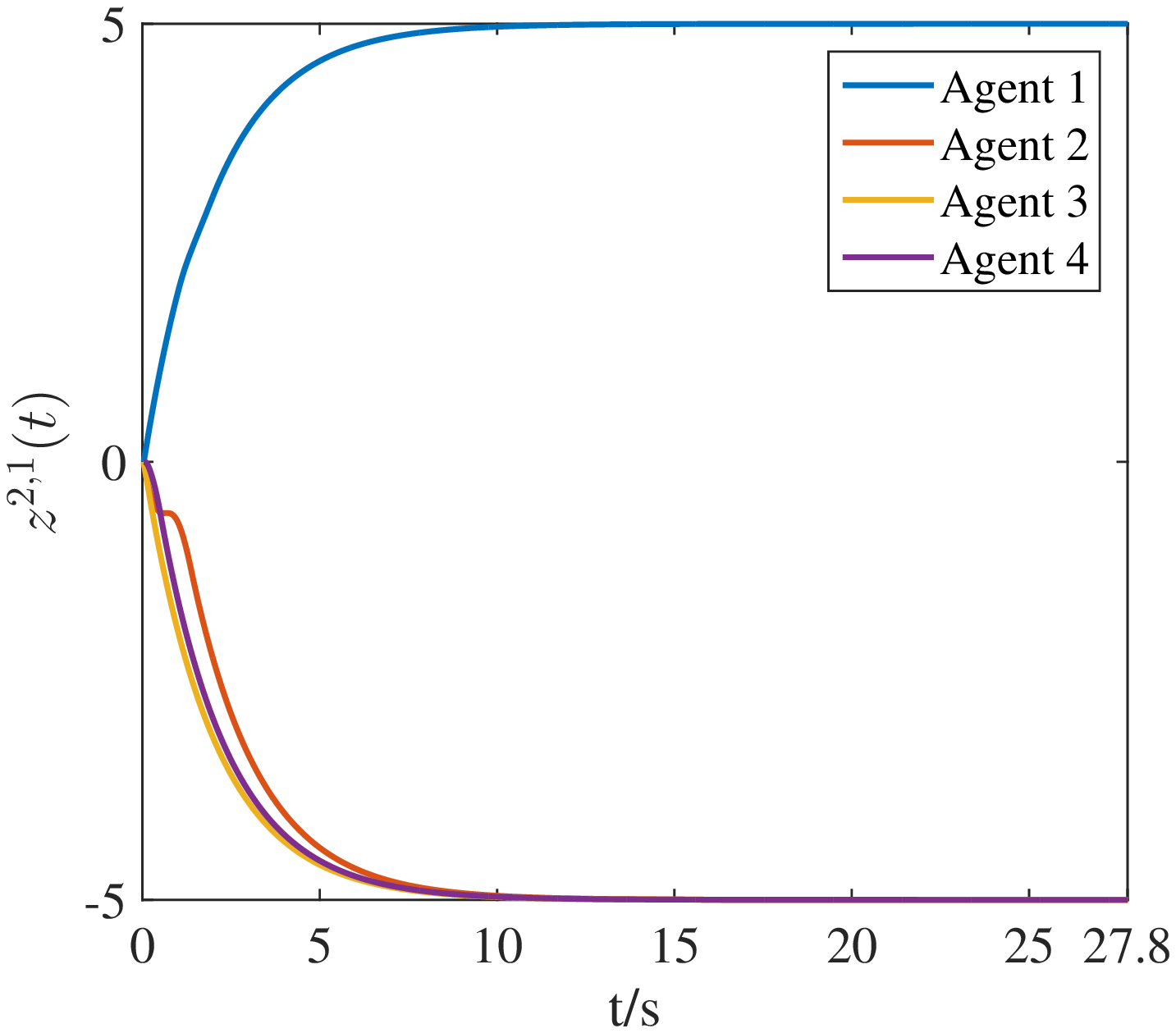}}
\subfigure{
\includegraphics[width=0.20\textwidth]{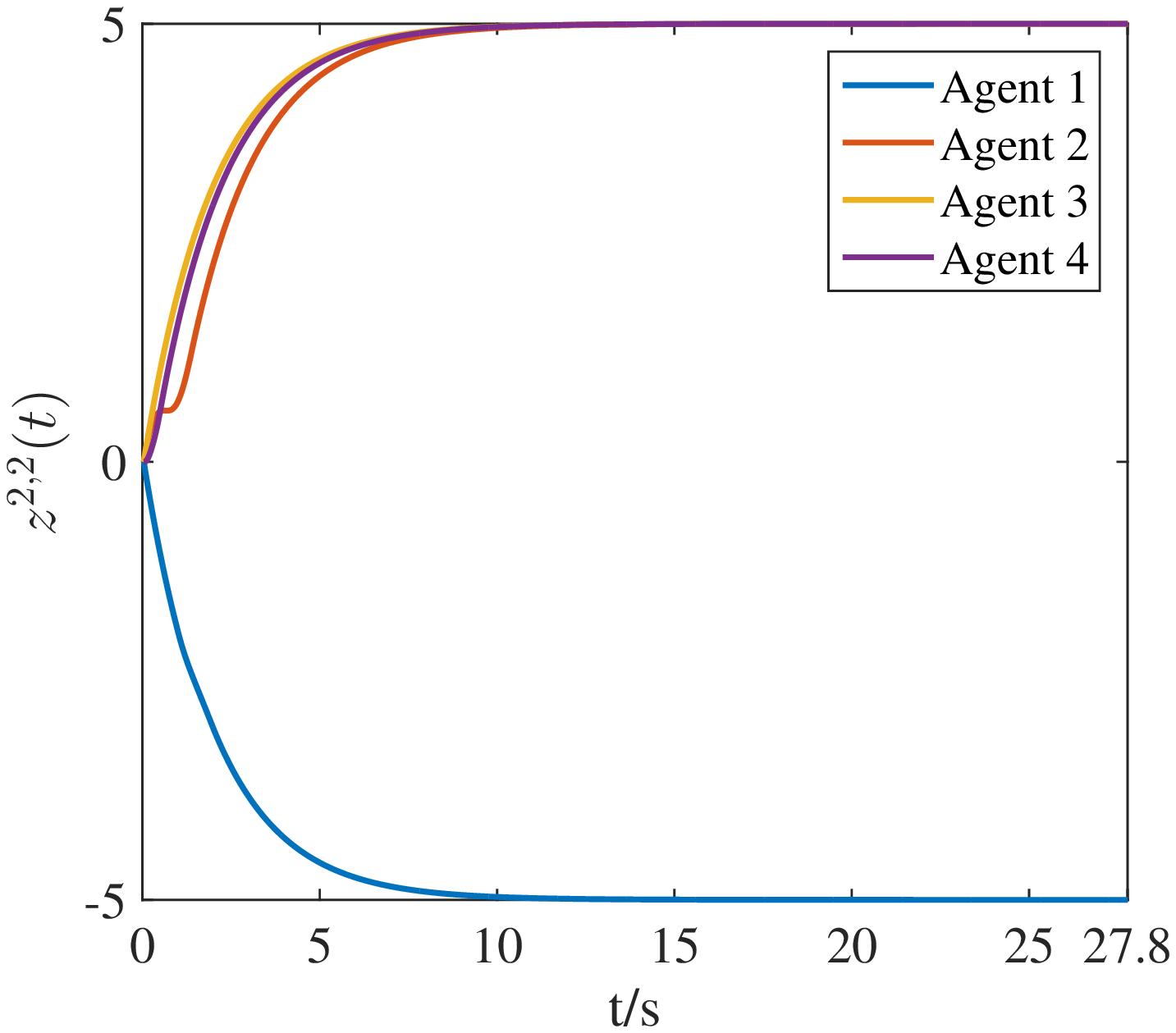}}
\caption{Trajectories of $z_{i}^{2}(t)$ for $i \in \lbrace 1, 2, 3, 4 \rbrace$ with algorithm ($\ref{Algorithm 2}$)}
\label{Fig.6}
\end{figure}

\begin{figure}
\centering
\subfigure{
\includegraphics[width=0.2\textwidth]{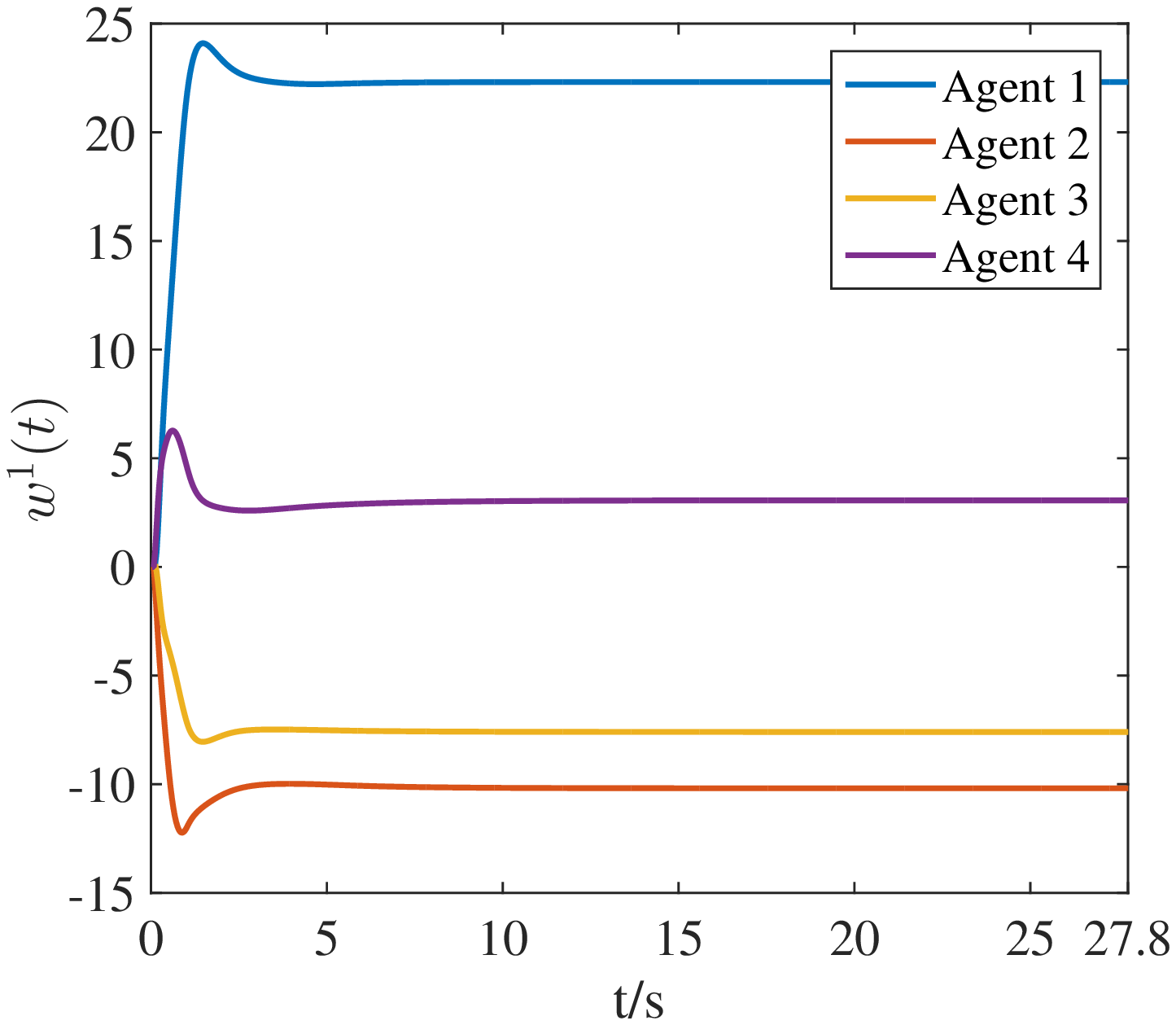}}
\subfigure{
\includegraphics[width=0.2\textwidth]{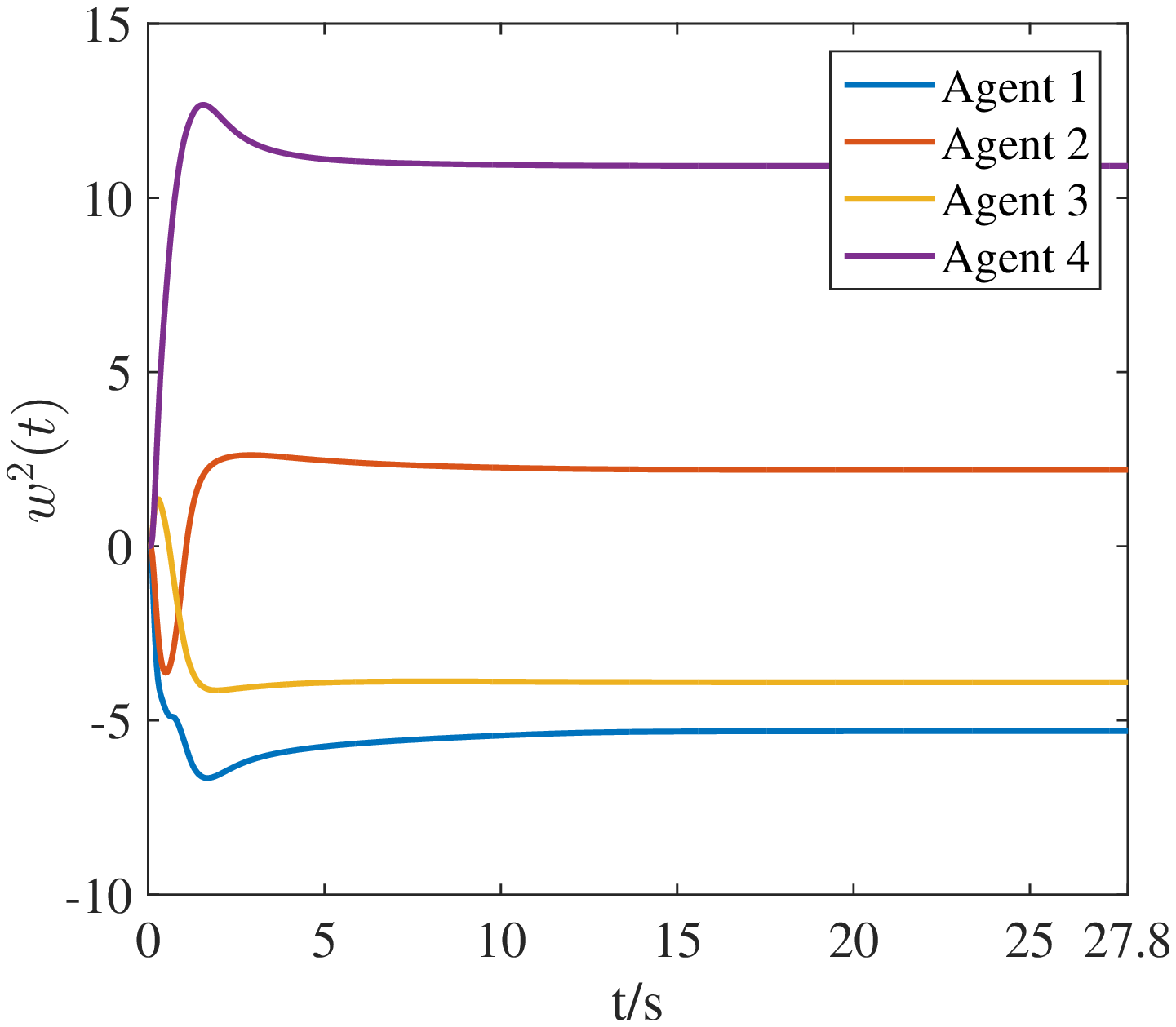}}
\caption{Trajectories of $w_{i}(t)$ for $i \in \lbrace 1, 2, 3, 4 \rbrace$ with algorithm ($\ref{Algorithm 2}$)}
\label{Fig.7}
\end{figure}

\section{Conclusion}
In this paper, a class of nonsmooth resource allocation problems with directed graphs was solved via two distributed multi-proximal operator based primal-dual algorithms. The second algorithm considered a distributed estimator of the left eigenvector $h$ corresponding to $\lambda_{1}(L_{nq}) = 0$. These two algorithms were smoothed thanks to the multi-proximal splitting. Moreover, the design of the second proposed algorithm can also give a new viewpoint to tackle many widely studied distributed constrained resource allocation problems. Future extensions will involve considering nonsmooth resource allocation problems with switching topologies and more complex communication situations such as time delay and packet losses.

\end{document}